\documentclass[12pt]{article}
\usepackage[final]{epsfig}
\usepackage{amsfonts, amsmath, amssymb, amsthm}
\usepackage{latexsym}
\usepackage{color}
\usepackage{multimedia}
\usepackage[T1]{fontenc}
\usepackage{tabls} 
\usepackage{multirow} 

\numberwithin{equation}{section}

\newtheorem{theorem}{Theorem}[section]

\newtheorem{example}{Example}[section]
\newtheorem{proposition}{Proposition}[section]
\newtheorem{definition}{Definition}[section]

\newcommand{\R}{\mathbb{R}} 

\newcommand{\Z}{\mathbb{Z}}

\newcommand{\be}{\begin{equation}}
\newcommand{\ee}{\end{equation}}
\newcommand{\calC}{{\mathcal{C}}}
\newcommand{\calP}{{\mathcal{P}}}

\newcommand{\calM}{{\mathcal{M}}}

\newcommand{\lambdabar}{{\overline{\lambda}}}
\newcommand{\rmF}{{\rm F}}
\newcommand{\rmK}{{\rm K}}
\newcommand{\rmT}{{\rm T}}
\newcommand{\rmC}{{\rm C}}

\begin{document}

\begin{center}
\begin{Large}
A simple counterexample to the Monge ansatz in multi-marginal optimal transport, 
convex geometry of the set of Kantorovich plans, and the Frenkel-Kontorova model \\
\end{Large}

\normalsize
\vspace{0.4in}

Gero Friesecke \\[1mm]
Faculty of Mathematics, Technische Universit\"at M\"unchen, {\tt gf@ma.tum.de} \\[2mm]
 
\end{center}

{\bf Abstract.} It is known from clever mathematical examples \cite{Ca10}
that the Monge ansatz may fail in continuous two-marginal optimal transport (alias optimal coupling alias optimal assignment) problems. Here we show that this effect already occurs for finite assignment problems with $N=3$ marginals, $\ell=3$ 'sites', and symmetric pairwise costs, with the values for $N$ and $\ell$ both being optimal. Our counterexample is a transparent consequence of the convex geometry of the set of symmetric Kantorovich plans for $N=\ell=3$, which -- as we show -- possess 22 extreme points, only 7 of which are Monge. These extreme points have a simple physical meaning as irreducible molecular packings, and the example corresponds to finding the minimum energy packing for Frenkel-Kontorova interactions. Our finite example naturally gives rise, by superposition, to a continuous one, where failure of the Monge ansatz manifests itself as nonattainment and formation of 'microstructure'. 
%
%
%

\section{Introduction} \label{sec:Intro} 
\begin{figure}
\begin{minipage}{0.49\textwidth}
\begin{center}
\includegraphics[width=0.9\textwidth]{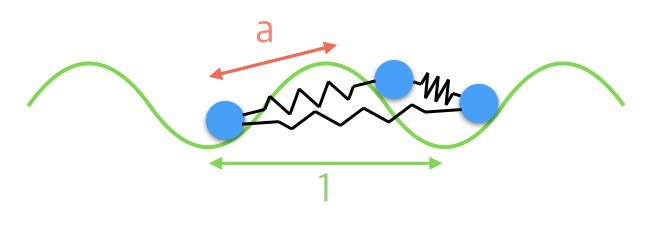} 
\end{center}
\end{minipage}
\begin{minipage}{0.49\textwidth}
\begin{center}
\vspace*{-0.2cm} 
\includegraphics[width=0.65\textwidth]{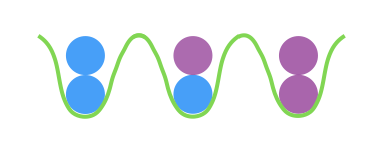} 
\end{center}
\end{minipage}
\vspace*{-5mm}
\caption{{\it Left:} The Frenkel-Kontorova model. Particles are linked by springs with positive equilibrium length $a$, and subject to an external potential of period $L$ (set here to $1$). The competition between springs and external potential leads to symmetry breaking, i.e. the preferred distance between particle 1 and 2 is different from that between particle 2 and 3. 
{\it Right:} Ground state of the total spring energy of an ensemble of $3$-particle molecules 
when the external potential is replaced by rigid confinement to three consecutive minima and we impose the constraint that each minimum must be occupied by the same number of particles. This is precisely the multi-marginal optimal transport problem from Example 1.1.}
\label{F:FKmodel}
\end{figure}

For which costs is the Monge ansatz justified in multi-marginal optimal transport? This deep question remains much less well understood than its two-marginal counterpart \cite{Vi09}. 
For interesting examples of Monge and non-Monge minimizers see, respectively, 
\cite{GS98, He02, Ca03, Pa11, CFK11, BDG12, CDD13, DGN15} and \cite{CN08, Pa10, FMPCK13, Pa13, CFP15, MP17, GKR18}. 
Our goal in this paper is to point out a fundamental difference between multi-marginal and two-marginal optimal transport. Namely, in the important case of finite state spaces and uniform marginals -- which arises by equi-mass discretization from any continuous optimal transport problem on $\R^d$ with absolutely continuous marginals \cite{CFM14} -- {\it the Monge ansatz is sufficient for 2 marginals, but not for 3 or more marginals.} The sufficiency for 2 marginals is the content of the remarkable Birkhoff-von Neumann theorem \cite{Bi46, vN53}. The insufficiency for 3 marginals is shown here; what is more, we construct by systematic arguments the lowest-dimensional counterexample possible for symmetric and pairwise costs, which turns out to occur for 3 marginals and 3 marginal states (or ``sites'', see below). The counterexample has a simple physical meaning related to the Frenkel-Kontorova model \cite{BK98} of solid state physics (see Figure \ref{F:FKmodel}), as we will of course explain. 

\begin{figure}
\begin{center}
\includegraphics[width=0.33\textwidth]{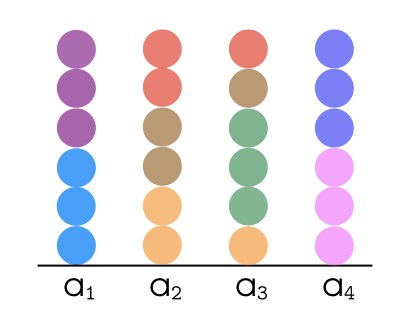}
\end{center}
\vspace*{-10mm}

\caption{A typical Kantorovich plan, visualized as a molecular packing. The picture corresponds to the plan $\gamma=\tfrac{1}{2}\delta_{a_2}\otimes\delta_{a_2}\otimes\delta_{a_3} + \tfrac{1}{3}\delta_{a_1}\otimes\delta_{a_1}\otimes\delta_{a_1} + \tfrac{1}{3}\delta_{a_4}\otimes\delta_{a_4}\otimes\delta_{a_4} + \tfrac{1}{6}\delta_{a_3}\otimes\delta_{a_3}\otimes\delta_{a_3}$. 
Since the probability coefficients are in the ratio $3:2:2:1$, we can visualize the whole plan as $3$ molecules with particle positions $a_2$, $a_2$, $a_3$, $2$ molecules with particle positions $a_1$, $a_1$, $a_1$, $2$ molecules with particle positions $a_4$, $a_4$, $a_4$, and one molecule all of whose particles are at $a_3$. The different colors serve to distinguish the different molecules.}
\label{F:1}
\end{figure}

For multi-index assignment problems, i.e. multi-marginal problems with finite state space, it has long been known that there exist ``non-integer vertices'' of the $-$ suitably renormalized $-$ polytope of Kantorovich plans \cite{Cs70, Kr07, LL14}. These can be shown \cite{Fi14} to correspond to non-Monge plans. By abstract duality principles this implies the existence of costs with unique non-Monge minimizers; but we are not aware of previous explicit examples, let alone ones of mathematically simple (e.g., pairwise and symmetric) and physically relevant form. 
\\[2mm]
{\bf Multi-marginal optimal transport.} 
Given a finite state space 
\be \label{space}
    X=\{a_1,...,a_\ell\}
\ee
consisting of $\ell$ distinct points, and a cost function $c\, : \, X^N\to\R\cup\{+\infty\}$, the $N$-marginal optimal transport problem in Kantorovich form consists of the following: Minimize the total cost
\be \label{Kant}
    \calC[\gamma] = \int_{X^N} c_N(x_1,...,x_N) \, d\gamma(x_1,...,x_N)
\ee
over probability measures $\gamma$ on $X^N$ with one-point marginals $\lambdabar$, i.e.
\be \label{Kantconstr} 
   \gamma(X^{i-1}\times A_i\times X^{N-i}) = \lambdabar(A_i) \mbox{ for all subsets }A_i\mbox{ of }X 
   \mbox{ and all }i=1,...,N, 
\ee
where $\lambdabar$ is the uniform measure on $X$, that is to say
\be \label{unif}
    \lambdabar = \sum_{i=1}^\ell \tfrac{1}{\ell} \, \delta_{a_i}.
\ee
Here and below, $\delta_{a_i}$ denotes the Dirac measure on the point $a_i$, and we use the common notation $\gamma\mapsto\lambdabar$ for the validity of eq.~\eqref{Kantconstr}. Probability measures $\gamma$ on $X^N$ satisfying \eqref{Kantconstr} are known as {\it Kantorovich plans}. 
\\[2mm]
The Monge form of the above optimal transport problem is to find $N$ permutations $\tau_1$, ..., $\tau_N \, : \, \{1,...,\ell\} \to \{1,...,\ell\}$ which minimize the cost
\be \label{Monge}
    \calC[\tau_1,...,\tau_N ] = \frac{1}{\ell} \sum_{\nu=1}^\ell c(a_{\tau_1 (\nu )},...,a_{\tau_N (\nu )}). 
\ee
This corresponds to making the special ansatz
\be \label{Mongeansatz}
    \gamma = \frac{1}{\ell} \sum_{\nu=1}^\ell 
    \delta_{a_{\tau_1(\nu)}}\otimes \cdots \otimes \delta_{a_{\tau_N(\nu)}} 
\ee
in the Kantorovich problem. Note that the requirement that the $\tau_i$ be permutations ensures the marginal condition \eqref{unif}. Kantorovich plans of this form are called {\it Monge states}. By re-ordering the sum in \eqref{Monge} and \eqref{Mongeansatz} one may assume $\tau_1=id$, but we prefer the above more symmetric formulation. 
\\[2mm]
The problems \eqref{Kant}--\eqref{unif} and \eqref{Monge}--\eqref{Mongeansatz} arise in many contexts including economics \cite{CE10, CMN10}, image processing \cite{AC11, RPDB12}, mathematical finance \cite{BHP13, GHT14}, optimal assignment problems (see the reviews \cite{Sp00, BDM12}), 
and -- more recently -- electronic structure \cite{CFK11, BDG12}. In the latter context, Kantorovich plans correspond to the joint probability distribution of electron positions in an $N$-electron molecule, $X$ is a collection of discretization points in $\R^3$, and the cost $c(x_1,...,x_N)$ is given by the Coulomb interaction energy $\sum_{1\le i<j\le N}|x_i-x_j|^{-1}$ between the electrons, with $|\cdot |$ being the euclidean distance.   
Partially inspired by this example, we propose a more general physical problem which is described by \eqref{Kant}--\eqref{unif}. We find that thinking about this physical problem provides novel and valuable intuition about \eqref{Kant}--\eqref{unif}. 
\\[2mm]
{\bf Molecular packing problem:} {\it Find the ground state of an ensemble of 
$N$-particle molecules confined to $\ell$ sites $a_1,...,a_\ell\in\R^d$ subject to the constraint that all sites must be occupied equally often. The molecules are composed of $N$ identical particles (``atoms''), and the cost $c$ which is to be minimized is the intramolecular interaction energy between the particles within a molecule.}
\\[2mm]
In this problem, the state of a single $N$-particle molecule is described by a Dirac measure $\delta_{x_1}\otimes \cdots \otimes \delta_{x_N}$, with $x_1,...,x_N$ denoting the positions of the $N$ particles, and the state of the ensemble is described by a superposition $\gamma=\sum_{\nu} p_\nu 
\delta_{x_1^{(\nu)}}\otimes \cdots \otimes \delta_{x_N^{(\nu)}}$ of single-molecule states, with the $p_\nu$ denoting occupation probabilities. The passage from single-molecule states to ensembles may be viewed as a physics analogue of the passage from pure to mixed states in game theory. A typical state of the ensemble can be visualized by filling the sites with as many $N$-atom molecules as needed to make the number of molecules in each state proportional to the state's occupation probability. See Figure \ref{F:1} for a typical state of such an ensemble, and Example 1.1 for an instructive example of a ground state. This visualization ignores the order in which the tensor factors $\delta_{x_i}$ appear, but uniquely characterizes the symmetrization $S\gamma$ of the state (with $S$ as defined in \eqref{S}). Note that the latter provides a physically appropriate description of an ensemble of molecules when the particles (``atoms'') within each molecules are indistinguishable. With the help of $S\gamma$ we can also give physical meaning to the marginal condition \eqref{Kantconstr}: for any probability measure $\gamma$ on $X^N$, 
\be \label{margphys}
   \mbox{$S\gamma$ satisfies \eqref{Kantconstr}} \;
   \Longleftrightarrow \; \begin{array}{c} \mbox{the associated molecular packing} \\
                                           \mbox{occupies all sites equally often}.
                          \end{array}
\ee

Due to the finiteness of the space $X$, probability measures on $X^N$ can be identified with order-$N$ tensors $(\gamma_{i_1...i_N})$ with components $\gamma_{i_1...i_N}\in\R$ and indices $i_1,...,i_N$ running from $1$ to $\ell$, via $\gamma_{i_1...i_N}:=\gamma(\{(a_{i_1},...,a_{i_N})\})$. It is instructive to write out our optimal transport problem in tensor notation in the special cases of two and three marginals. 

{\it Two-marginal case:} For $N=2$, the Kantorovich optimal transport problem corresponds via $c_{ij}:=c(a_i,a_j)$ and $\mu_{ij}:=\ell\gamma(\{(a_i,a_j)\})$ to the celebrated Birkhoff-von Neumann problem: for given cost coefficients $c_{ij}$, 
\be \label{BvN}
  \mbox{Minimize }\sum_{i,j=1}^\ell c_{ij}\mu_{ij}
\ee
over the $\mu_{ij}$ subject to the constraint that $(\mu_{ij})_{i,j=1}^\ell$ is a bistochastic matrix, i.e. that
\be \label{BvNconstr}
  \mu_{ij}\ge 0 \mbox{ for all }i,j, \;\; \sum_{i}\mu_{ij}=1 \mbox{ for all }j, \;\; \sum_{j}\mu_{ij}=1 \mbox{ for all }i.
\ee
(The prefactor $\ell$ in the definition of the $\mu_{ij}$ just serves to adapt the probability coefficients $\gamma(\{(a_i,a_j)\})$, which sum to $1$, to the normalization convention in economics and game theory, where one instead wants the partial sums in \eqref{BvNconstr} to be $1$.) 
The Monge problem corresponds to the sparse ansatz that $(\mu_{ij})$ is a permutation matrix, that is to say each row and each column contains exactly one $1$ and $\ell\! - \! 1$ zeros.

{\it Three-marginal case:} For $N=3$, the Kantorovich optimal transport problem corresponds via $c_{ijk}:=c(a_i,a_j,a_k)$ and $\mu_{ijk} = \ell \gamma(\{(a_i,a_j,a_k)\})$ to the axial three-index assignment problem: 
\be \label{A3IAP}
  \mbox{Minimize }\sum_{i,j,k=1}^\ell c_{ijk}\mu_{ijk} 
\ee
subject to the constraint that $(\mu_{ijk})$ is a tristochastic $\ell\times\ell\times\ell$ tensor of order three, i.e.
\be \label{A3IAPconstr}
  \begin{array}{ll} \mu_{ijk}\ge 0 \mbox{ for all }i,j,k, 
                  & \sum_{i,j} \mu_{ijk}=1 \mbox{ for all }k, \\ 
                    \sum_{i,k} \mu_{ijk}=1 \mbox{ for all }j=1, 
                  & \sum_{j,k} \mu_{ijk}=1 \mbox{ for all }i=1. 
  \end{array}
\ee
The Monge problem corresponds to the sparse ansatz that each of the 3$\ell$ ``planes'' associated with the above sum constraints, i.e. $(\mu_{ijk})_{i,j=1}^\ell$, $(\mu_{ijk})_{i,k=1}^\ell$, $(\mu_{ijk})_{j,k=1}^\ell$, contains exactly one $1$ and $\ell^2\! - \! 1$ zeros. 
\\[2mm]
{\bf The counterexample.} 
The Birkhoff-von Neumann theorem \cite{Bi46, vN53} says that for $N=2$, the extreme points of the polytope of bistochastic $\ell\times\ell$ matrices are precisely the permutation matrices. In particular, for any cost function $c \, : \, X\times X\to\R$ the Kantorovich problem \eqref{Kant}--\eqref{unif} possesses a minimizer which is of Monge form \eqref{Monge}. 
\\[2mm]
For $N=3$ and $\ell=3$ we will show that the polytope of tristochastic $\ell\times\ell\times\ell$ tensors of order three possesses extreme points which are not of Monge form; in particular, there exist cost functions $c \, : \, X\times X\times X\to\R$ such that none of the minimizers of the Kantorovich problem \eqref{Kant}--\eqref{unif} is of Monge form \eqref{Monge}. 
As we will see, this phenomenon even occurs in the class of {\it pairwise and symmetric} costs, i.e.
\be \label{pairwsymm}
   c(x,y,z) = v(x,y) + v(x,z) + v(y,z) 
   \mbox{ for some }v\, : \, X\times X \to \R \mbox{ with } v(x,y)=v(y,x). 
\ee  
%
%
\begin{example} \label{E:counter} (Optimal packing of Frenkel-Kontorova molecules) Consider the Kantorovich problem \eqref{Kant}--\eqref{unif} with $N=\ell=3$ and $X$ given by three equi-spaced points on the real line, i.e. $X=\{1,2,3\}\subset\R$ (physically: consider an ensemble of $3$-particle molecules confined to the sites $1$, $2$, $3\in\R$ subject to the constraint that all sites must be occupied equally often). For the  cost
\be \label{counterex}
   C[\gamma] = \int_{X^3}\bigl( v(|x-y|) + v(|y-z|) + v(|x-z|)\bigr)\, d\gamma(x,y,z), \;\; v(r)=(r - a)^2, \;\; a=\tfrac{3}{4}
\ee
(physically: when the particles within a molecule are mutually connected by springs of equilibrium length $a=\tfrac{3}{4}$), the unique minimizer is given by $\gamma_* = S \gamma$ where $\gamma=\tfrac{1}{2}(\delta_1\otimes\delta_1\otimes\delta_2 + \delta_2\otimes\delta_3\otimes\delta_3)$ and $S$ is the symmetrization operator \eqref{S}. This $\gamma_*$ is not a Monge state. It is not a symmetrized Monge state either. 
\end{example}

Example 1.1 is the ``simplest example possible'', in the following precise sense: $N$ and $\ell$ are minimal (see Remark 2 following Theorem \ref{T:1}); and $v$ cannot be taken to be monotone (see Examples 4.3 and 4.4 in section \ref{S:gallery}). 
\\[2mm]
Physically, the above minimizer is very intuitive, and can be heuristically derived as follows.  
In its original form, the Frenkel-Kontorova (FK) model \cite{BK98} describes the equilibria and excitations of a 1D chain of particles (``atoms'') linked by springs of positive equilibrium length and subject to a periodic external potential (describing the interaction of the chain with a ``substrate''); in our case the chain is a ring, with each atom linked to the two others. See Figure \ref{F:FKmodel}, left panel. Now Example 1.1 corresponds to the limiting situation when particle positions in the ``holes'' (minima of the periodic external potential) are not just energetically favoured but rigidly enforced. Moreover we require all holes to be filled equally, and the only quantity that remains to be minimized is the interaction energy due to the springs. The FK model, and its adaptation here, favours symmetry breaking (see again Figure \ref{F:FKmodel}, left panel) when the spring equilibrium length is a bit smaller than the distance between the holes. Hence one expects that the FK ground state, when enforcing particle positions in the ``holes'', consists of two particles occupying the same hole and the third particle a neighbouring hole (see Section \ref{S:geo} for a mathematical proof of this fact). But it is intuitively obvious that the -- up to repetition only -- way to fill $3$ neighbouring holes uniformly by ground state molecules is the packing in Figure \ref{F:FKmodel}, right panel. Next we translate this packing into a Kantorovich plan according to the prescription in Figure \ref{F:1}. The blue state corresponds to $\delta_{1}\otimes\delta_1\otimes\delta_2$, since the first two blue particles are in hole $1$ and the third blue particle is in hole $2$. Likewise, the purple state corresponds to $\delta_{2}\otimes\delta_3\otimes\delta_3$. The whole ensemble corresponds 
to $\gamma=\tfrac{1}{2}\delta_{1}\otimes\delta_1\otimes\delta_2 + \tfrac{1}{2}\delta_{2}\otimes\delta_3\otimes\delta_3$, where the prefactors $\tfrac{1}{2}$ give the probability of the ensemble to be in either state. Finally, since the three holes are occupied equally often, the symmetrization $\gamma_*:=S\gamma$ satisfies the marginal constraint \eqref{Kantconstr} (see \eqref{margphys}). 

The key to deriving the ground state mathematically -- and to understand how one might come up with Example 1.1 in the first place --  
is to avoid intransparent ad hoc calculations and instead study the convex geometry of the set of Kantorovich plans (see Section \ref{S:geo}).
\\[2mm]
The overall significance of Example 1.1 is that it destroys -- even for pairwise, symmetric, simple costs -- the hope that the low-dimensional Monge ansatz \eqref{Mongeansatz} can be used to numerically tackle the high-dimensional Kantorovich problem \eqref{Kant}--\eqref{Kantconstr} for large $N$ and $\ell$, as would be desirable e.g. in applications to electronic structure; not that the number of unknowns in \eqref{Kant}--\eqref{Kantconstr} grows combinatorially in both $N$ and $\ell$. An almost as low-dimensional ansatz which cures the insufficiency of Monge states at the expense of only $2\ell$ additional unknowns independently of $N$ is presented in \cite{FV18}, where Example 1.1 was announced. 
\section{Convex geometry of the set of Kantorovich plans} \label{S:geo}
First, following \cite{FMPCK13} we introduce some relevant sets of probability measures. This may look technical at first, but it makes optimal transport {\it simpler}: the only thing one really needs to understand is the convex geometry of these (universal, i.e. cost-function-independent) sets. Let $X$ be the finite state space \eqref{space}, and let $\calP(X^N)$ denote the set of probability measures on $X^N$. The set of admissible competitors in the Kantorovich problem \eqref{Kant}--\eqref{unif} are the probability measures with uniform one-point marginal,
\be \label{Plambdabar}
   \calP_{\lambdabar}(X^N)  := \{ \gamma \in \calP(X^N) \, : \, \gamma\mapsto\lambdabar \}.
\ee
Due to finiteness of the state space $X$, $\calP_{\lambdabar}(X^N)$ is a convex polytope, i.e. a compact convex finite-dimensional set with only a finite number of extreme points. Recall that a point $x$ in a convex set $K$ is an extreme point if, whenever $x=\alpha x_1 + (1-\alpha)x_2$ for some $x_1$, $x_2\in K$ and some $\alpha\in(0,1)$, we have that $x_1=x_2=x$; recall also that, by Minkowski's theorem, any compact convex finite-dimensional set is the convex hull of its extreme points. 

\begin{figure}
\begin{center} 
\includegraphics[width=0.75\textwidth]{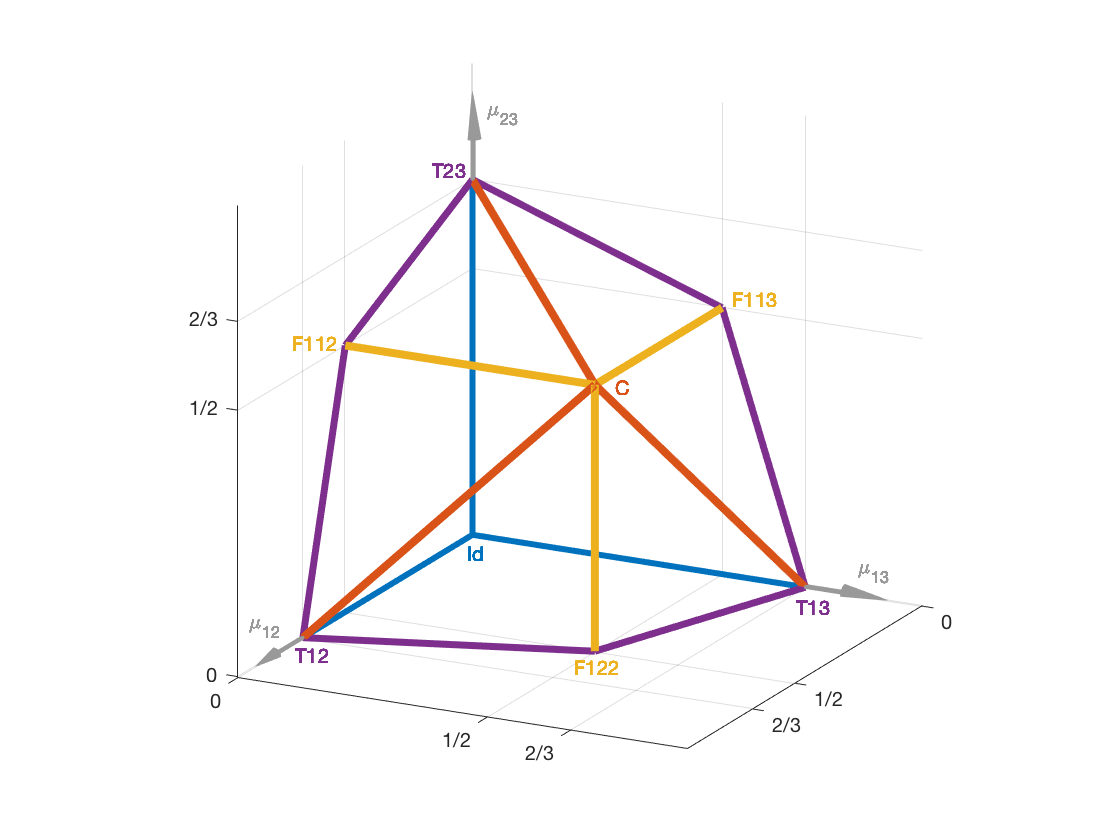}
\end{center}
\vspace*{-12mm}

\caption{The reduced Kantorovich polytope $\calP_{N-rep,\lambdabar}(X^2)$ for three marginals and three states ($N=\ell=3$). The coordinate axes indicate the upper triangular matrix elements in the representation \eqref{matrep}. These matrix elements uniquely specify any element of the Kantorovich polytope, due to the symmetry and bistochasticity of the matrix $(\mu_{ij})$. The red, blue and purple vertices correspond to symmetrized Monge states but the yellow vertices do not. The nomenclature for the vertices is explained in Remark 3.}  
\label{F:KantPoly}
\end{figure}

\begin{figure}
\begin{center}
\vspace*{-1cm}
\includegraphics[width=0.6\textwidth]{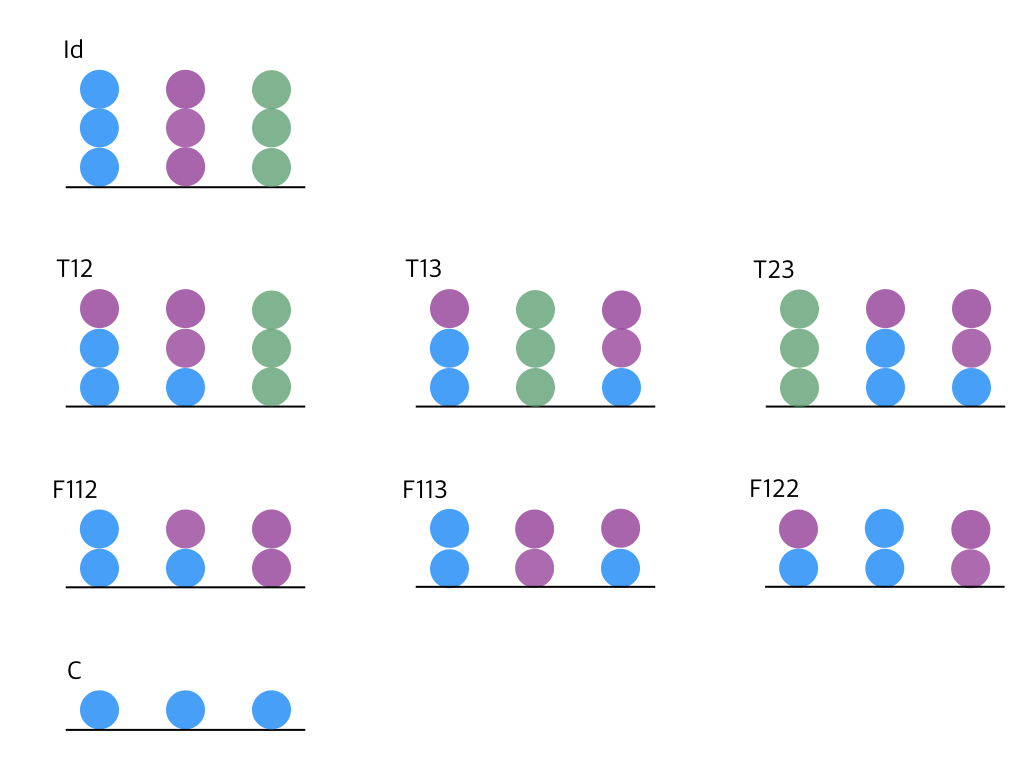} \\
\vspace*{-2.5mm}
\includegraphics[width=0.6\textwidth]{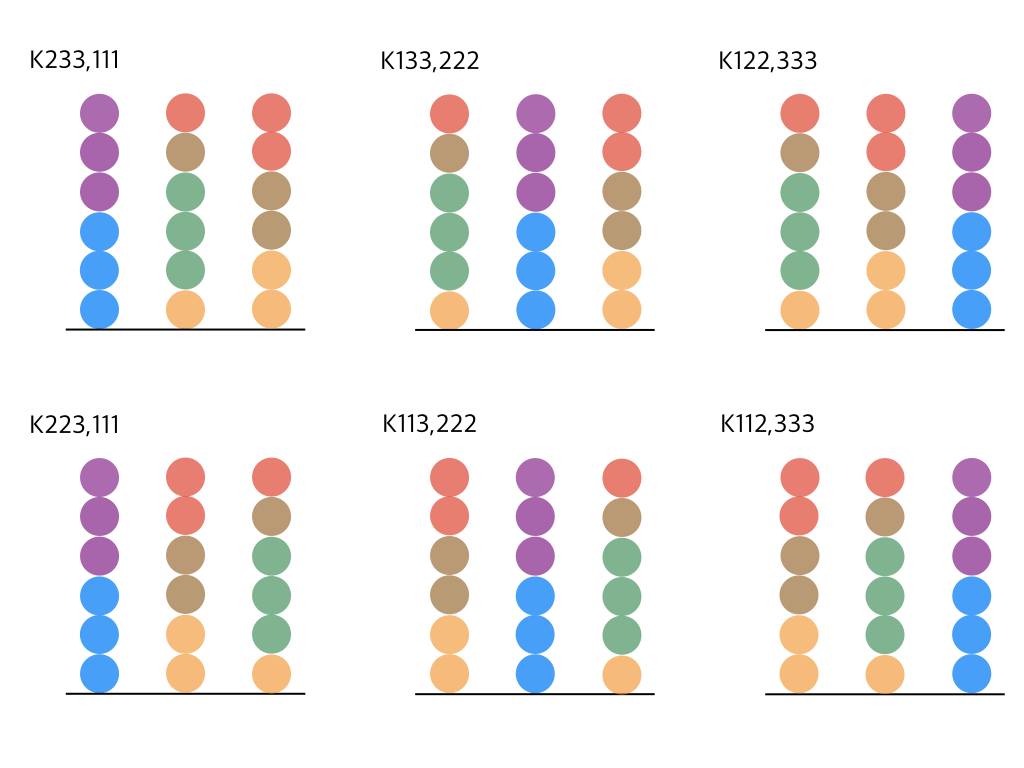} \\
\vspace*{-7.5mm} 
\includegraphics[width=0.6\textwidth]{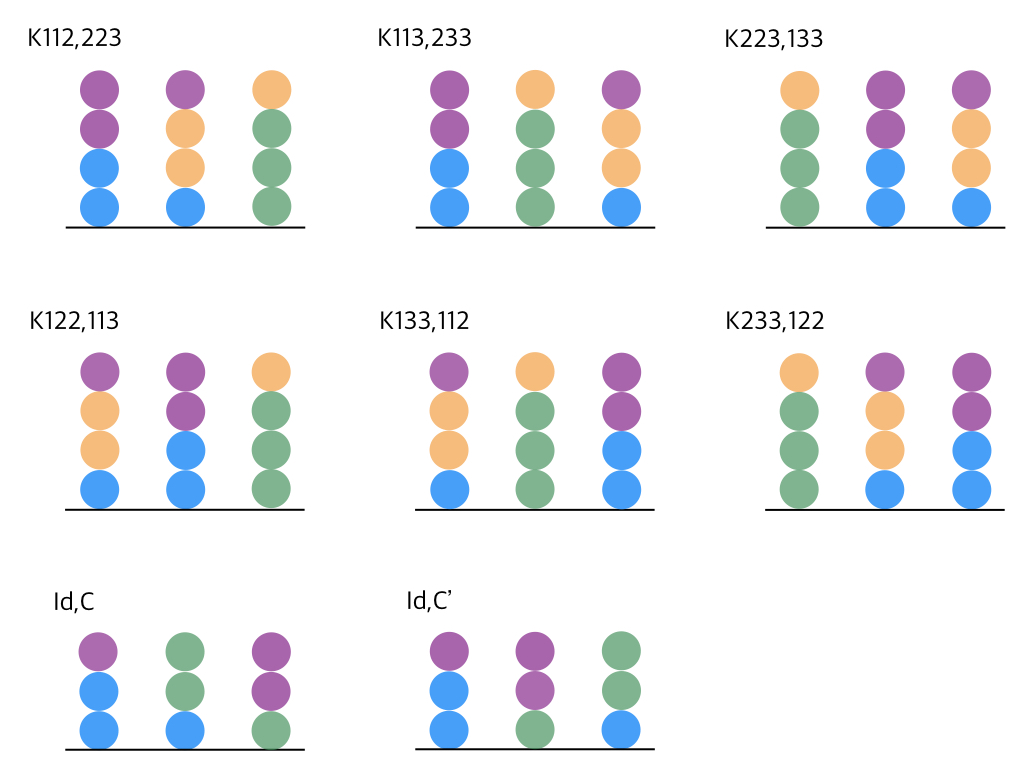}
\end{center}
\vspace*{-2.5mm} 
\caption{Physical meaning of the 22 extreme points of the symmetric Kantorovich polytope $\calP_{sym,\lambdabar}(X^3)$ for 3 marginals and 3 sites ($N=\ell=3$) as molecular packings. Each  Kantorovich plan corresponds to a uniform filling of the $3$ sites by $3$-particle molecules; extremal plans correspond to irreducible fillings (see \eqref{irred}). See Figure \ref{F:1} for how to translate the pictures into the plans in Table \ref{Table:1}.}
\label{F:KantStatesIconized}
\end{figure}

A simplification arises because for physical reasons, in this paper we are only interested in symmetric and pairwise costs, i.e. costs $c\, : \, X^N\to\R\cup\{+\infty\}$ such that
\be \label{pairwsymm'} 
   c(x_1,...,x_N) = \sum_{1\le i<j\le N} v(x_i,x_j) \mbox{ for some }v\, : \, X\times X\to\R\cup\{+\infty\} \mbox{ with }v(x,y)=v(y,x).
\ee
For such costs we trivially have that
\be \label{costid}
   \calC[\gamma] \; = \; \calC[S\gamma] 
   \; = \; \underbrace{\mbox{${N\choose 2}$} \int_{X\times X} v(x,y) \, d\bigl(M_2S\gamma\bigr)(x,y)}_{=:\calC'[M_2 S \, \gamma ]}
\ee
where $S$ is the symmetrizer and $M_2$ is the map from $N$-point probability measures to their two-point marginals, i.e.
\be \label{S}
   (S\gamma )(A_1 \times\cdots\times A _N) = \frac{1}{N!} \sum_{\sigma\in S_N} 
   \gamma\bigl( A_{\sigma(1)} \times \cdots \times A_{\sigma(N)}\bigr) 
   \mbox{ for all }A_1,...,A_N\subseteq X,
\ee
where $S_N$ denotes the group of permutations $\sigma \, : \{ 1,...,N\} \to \{ 1,...,N\}$, and
\be \label{M2} 
   (M_2\gamma)(A) = \gamma(A\times X^{N-2}) \mbox{ for all }A\subseteq X^{N-2}. 
\ee
The identity \eqref{costid} says that the cost $\calC[\gamma]$ depends on $\gamma$ only through its symmetrization $S\gamma$, and moreover only on the two-point marginal of the latter. This allows the following simple but fruitful reformulations of the Kantorovich problem \eqref{Kant}--\eqref{unif} as optimization problems over successively lower-dimensional polytopes obtained by successively applying $S$ and $M_2$, $S\calP_\lambdabar(X^N)=:\calP_{sym,\lambdabar}(X^N) $ and $M_2\calP_{sym,\lambdabar}(X^N)=:\calP_{N-rep,\lambdabar}(X^2) $. Before stating the reduced problems, let us comment on these two sets and explain the above notation. 
\begin{definition} A probability measure $\gamma\in\calP(X^N)$ is called symmetric if and only if $\gamma(A_1\times ... \times A_N) = \gamma(A_{\sigma(1)}\times ... \times A_{\sigma(N)})$ for all subsets $A_i$ of $X$ and all permutations $\sigma$.
\end{definition}
Clearly, a probability measure $\gamma\in\calP(X^N)$ is symmetric if and only if $\gamma=S\gamma$.  Hence the first of our two lower-dimensional sets is the set of symmetric $N$-point probability measures with uniform marginal,
\be \label{newset1}
  \calP_{sym,\lambdabar}(X^N) = \{ \gamma\in\calP(X^N) \, : \, \gamma \mbox{ symmetric, }
  \gamma\mapsto\lambdabar\}.
\ee
\begin{definition} {\rm \cite{FMPCK13}} A probability measure $p$ on $X^2$ is called {\it $N$-representable} (for some $N\ge 2$) if there exists a symmetric probability measure $\gamma$ on $X^N$ such that $M_2\gamma = p$. 
\end{definition}
Hence the second set is the set of $N$-representable two-point probability measures with uniform marginal, 
\be \label{newset2}
  \calP_{N-rep,\lambdabar}(X^2) = \{ p\in\calP(X^2) \, : \, p \mbox{ $N$-representable, }p\mapsto\lambdabar\}.
\ee
It follows from \eqref{costid} that
\be \label{minid}
   \min_{\gamma\in\calP_{\lambdabar}(X^N)} \calC[\gamma] \; = \; 
   \min_{\gamma\in\calP_{sym,\lambdabar}(X^N)} \calC[\gamma] \; = \; 
   \min_{p\in\calP_{N-rep,\lambdabar}(X^2)} \calC'[p].    
\ee
Thus symmetric costs are ``dual'' to the smaller polytope of symmetric $N$-point probability measures
with uniform marginal; symmetric and pairwise costs are ``dual'' to the even smaller polytope of $N$-representable two-point probability measures with uniform marginal. 
\\[2mm]
In the following, we call the sets $\calP_{sym,\lambdabar}(X^N)$ (see \eqref{newset1}) and $\calP_{N-rep,\lambdabar}(X^2)$ (see \eqref{newset2}), respectively, the {\it symmetric Kantorovich polytope for N marginals and $\ell$ states} and the {\it reduced Kantorovich polytope for $N$ marginals and $\ell$ states}.
\\[2mm]
As already mentioned in the introduction, we may represent any probability measure $p$ on $X^2$ with uniform marginals via
\be \label{matrep}
   \mu_{ij} = \ell \, p(\{(a_i,a_j)\} ) \;\;\; (i,j=1,...,\ell)
\ee
by a bistochastic $\ell\times\ell$ matrix $(\mu_{ij})_{i,j=1}^\ell$. In this matrix representation, the reduced Kantorovich polytope $\calP_{N-rep,\lambdabar}(X^2)$ becomes a subset of the Birkhoff polytope of all bistochastic $\ell\times\ell$ matrices. 
\\[2mm]
We will frequently use the shorthand notation 
\be \label{sh}
    \delta_i := \delta_{a_i}, \;\;\; \delta_{i_1...i_N} := \delta_{a_{i_1}}\otimes\cdots\otimes\delta_{a_{i_N}}.
\ee
\begin{theorem} \label{T:1} The reduced Kantorovich polytope $\calP_{N-rep,\lambdabar}(X^2)$ for $3$ marginals and $3$ states ($N=\ell=3$) has precisely $8$ extreme points, given by $M_2\gamma$ where $\gamma$ is one of the following Kantorovich plans from Table \ref{Table:1}:
\be \label{corners}
    {\rm Id, \; T12, \; T13, \; T23, \; C, \; F112, \; F113, \; F122.}
\ee
(See Figure \ref{F:KantPoly} for a picture of the polytope, Figure \ref{F:KantStatesIconized} for the physical meaning of the 8 plans, and Remark 3 for an explanation of nomenclature.) Moreover for each extreme point $p=M_2\gamma$, $\gamma$ is the unique element of the symmetric Kantorovich polytope $\calP_{sym,\lambdabar}(X^N)$ such that $p=M_2\gamma$. The first five of the plans \eqref{corners} are symmetrized Monge plans (i.e. of the form $S\tilde{\gamma}$ for some Monge state $\tilde{\gamma}$) but the last three are not.
\end{theorem}

Some remarks are in order.
\\[2mm]
1) The most important message of Figure \ref{F:KantPoly} is that some but not all of the vertices are Monge. Hence there exist costs for which the Monge ansatz is wrong. Moreover such costs can even be found in the symmetric pairwise class \eqref{pairwsymm'}, with which the pictured polytope is in duality (see \eqref{minid}).
\\[2mm]
2) \label{R:1} The numbers $N$ and $\ell$ in Theorem \ref{T:1} are minimal for the existence of non-Monge vertices. For $N=2$ such vertices are ruled out by the Birkhoff-von Neumann theorem. And for $N=3$, $\ell=2$ it is easy to show using the results in \cite{FMPCK13} that the Kantorovich polytope is the line segment with the two endpoints $M_2\gamma_1$, $M_2\gamma_2$, where
$\gamma_1= \tfrac{1}{2}(\delta_{111}+\delta_{222})$ and $\gamma_2 = S \tfrac{1}{2}(\delta_{112} + \delta_{221})$; but these $\gamma$'s are symmetrized Monge states. 
\\[2mm]
3) The above picture provides fruitful and hitherto lacking intuition about {\it when} the Monge ansatz is correct. Costs favouring small off-diagonal entries will drive the system to the bottom left back corner, which is Monge; this is suggestive of the fundamental continuum result of Gangbo and \'{S}wi\k{e}ch that the Monge ansatz is correct for the multi-marginal Wasserstein cost $c(x_1,...,x_N)=\sum_{i<j}|x_i-x_j|^2$ on the Euclidean space $X=\R^d$. Costs favouring large diagonal entries will drive the system to the top right front corner, which is also Monge; this provides a hint towards the important continuum result of Colombo, DiMarino, and DePascale \cite{CDD13} which justifies the Monge ansatz for the multi-marginal 1D Coulomb cost $\sum_{i<j}|x_i-x_j|^{-1}$. ``Intermediate'' costs favouring some off diagonal elements to be large and others to be small should favour vertices such as $\rmF 112$ which are seen in the picture to minimize $\mu_{13}$ and also maximize the sum $\mu_{12}+\mu_{23}$; in particular, this suggests that costs favouring symmetry breaking (like the Frenkel-Kontorova cost in Figure \ref{F:FKmodel}) may lead to non-Monge vertices, and led to the design of Example 1.1.
For a systematic translation of these insights into rigorous results for interesting classes of costs see Section \ref{S:gallery}. 
\\[2mm]
3) \label{R:2} The nomenclature in Theorem \ref{T:1}, Figure \ref{F:KantPoly} and Table \ref{Table:1} is as follows. Monge states are labelled by the underlying permutations $\tau\, : \, \{1,2,3\}\to\{1,2,3\}$, as follows. ${\rm Id}$, ${\rm T}_{ij}$ ($i<k$), ${\rm C}$, and ${\rm C'}$ stand, respectively, for the identity, the transposition of elements $i$ and $j$, the cyclic permutation $1\mapsto 2$, $2\mapsto 3$, $3\mapsto 1$, and its inverse. Recall that, for any Monge state \eqref{Monge}, by re-ordering the sum we can take $\tau_1$ to be the identity. The Monge state associated with the three permutations ${\rm Id}$, $\tau$, $\tau^2$ is denoted by $\tau$, as is its symmetrization. Thus, for example, ${\rm T}12$ stands for the state 
$S\tfrac{1}{3}\sum_{\nu=1}^3
\delta_{Id(\nu)}\otimes\delta_{{\rm T}_{12}(\nu)}\otimes\delta_{{\rm T}_{12}^2(\nu)}$. 
Interestingly, such multi-marginal Monge states generated by a single permutation have previously appeared in continuous optimal transport problems \cite{GM13, CDD13}. Any Monge state not generated by a single permutation, i.e. associated with three permutations ${\rm Id}$, $\tau$, $\tau'$ with $\tau'\neq \tau^2$, is denoted by $\tau,\tau'$, as is its symmetrization. We call the fundamental Kantorovich plans of non-Monge form which appear in Theorem \ref{T:1} and Figure \ref{F:KantPoly} Frenkel-Kontorova states, on account of their occurrence as minimizers in the Frenkel-Kontorova problem in Example 1.1; hence the letter $\rmF$. The subsequent digits indicate the first Dirac mass appearing; note that the second Dirac mass is, up to the order of appearance of the indices, uniquely determined by the marginal condition, and hence unique under the convention that indices appear in nondecreasing order. Thus, for example, $\rmF 112$ stands for the ground state $\frac{1}{2}S(\delta_{112} + \delta_{233})$ in Example 1.1. Finally, the more complicated  Kantorovich plans of non-Monge form appearing in Table \ref{Table:1} are denoted with the letter $\rmK$ followed by the index sequences of the biggest two Dirac masses; under our convention that indices apper in nondecreasing order, the remaining Dirac mass is again unique. Thus, e.g., $\rmK 233,111$ stands for the state $\frac{1}{2}S\delta_{233} + \tfrac{1}{3}\delta_{111} + \tfrac{1}{6}\delta_{222}$.

\begin{figure}
\begin{center}
\vspace*{-1cm}
\includegraphics[width=0.6\textwidth]{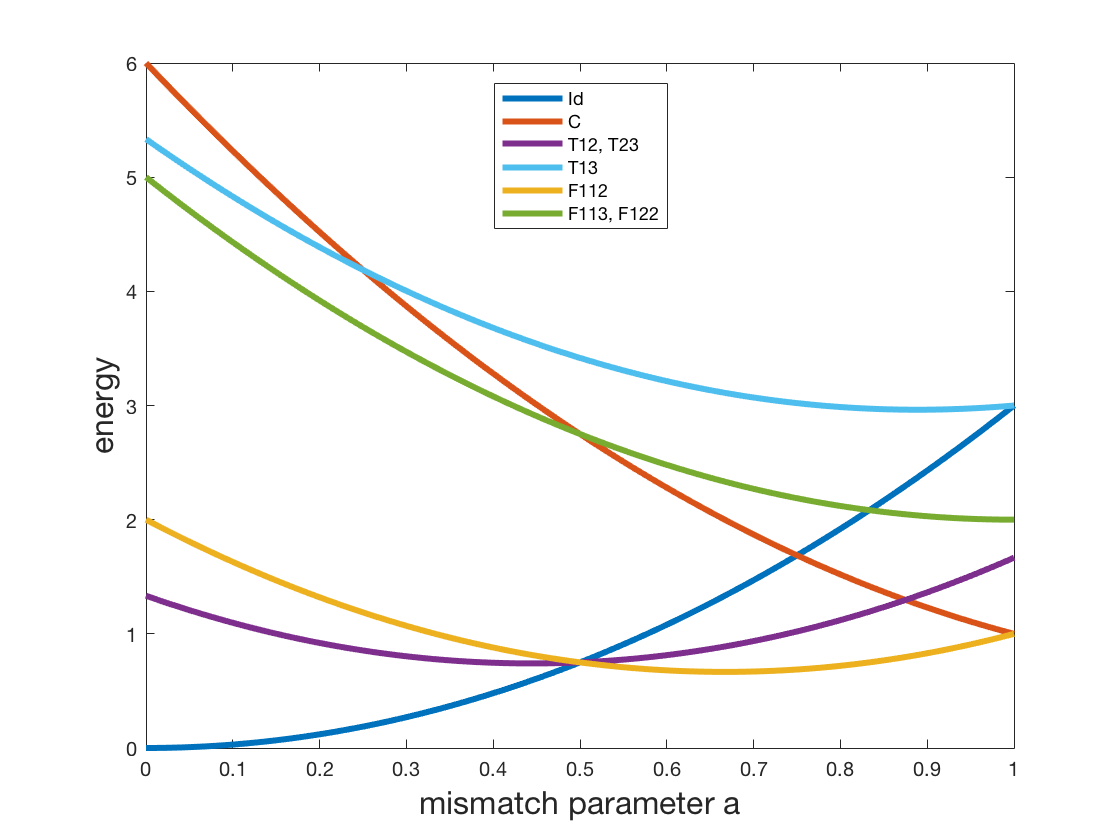}
\end{center}
\vspace*{-2.5mm} 
\caption{Frenkel-Kontorova energy \eqref{counterex} of the eight Kantorovich plans with extremal $2$-point marginal from Figure \ref{F:KantPoly}, as a function of the mismatch parameter $a$ (see Figure \ref{F:FKmodel}). As $a$ crosses $\tfrac{1}{2}$, the minimizer changes from the Monge plan  ${\rm Id}$ to the non-Monge plan $\rmF 112$.}
\label{F:FKresults}
\end{figure}

\begin{figure}
\begin{center}
\includegraphics[width=0.75\textwidth]{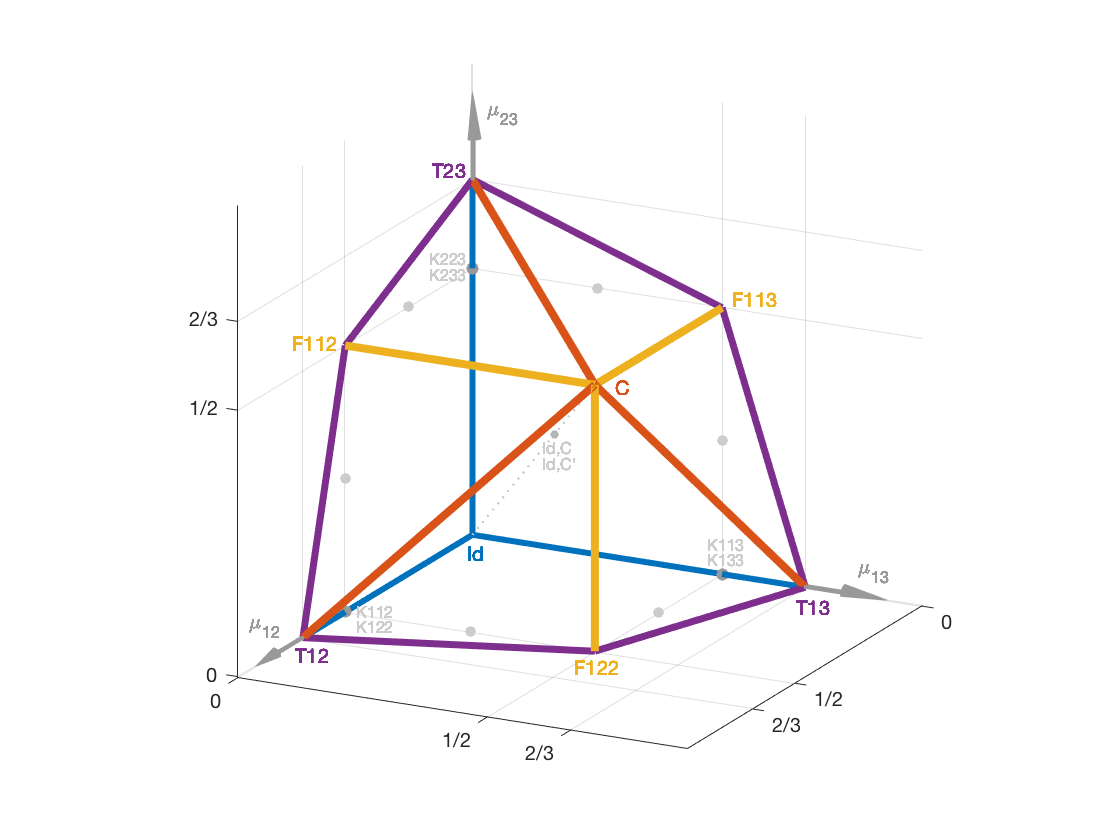} \\
\end{center}
\vspace*{-2.5mm} 
\caption{The 2-point marginals of the 22 extreme points of the symmetric Kantorovich polytope $\calP_{sym,\lambdabar}(X^3)$ for 3 marginals and 3 sites ($N\! =\!\ell\!=\! 3$). Coordinate axes are as in Fig.~\ref{F:KantPoly}. Here $Kiij$ stands for the state $Kiij,kkk$, and the unannotated grey points correspond to the states $Kiij,jjk$.}
\label{F:KantPoly2}
\end{figure}

\begin{table}[htbp]
\hspace{5mm}
\resizebox{0.9\textwidth}{!}{
  \begin{tabular}{|l|l|l|l|l|}
 \hline 
Name & Kantorovich plan $\gamma$ & Monge? & 2-point marginal times $\ell$ & $\begin{array}{l}\mbox{2-point marg.} \\ \mbox{extremal?}\end{array}$  \\
    \hline \hline
    Id & $\tfrac{1}{3}(\delta_{111}+\delta_{222}+\delta_{333})$ & Yes & $\delta_{11}+\delta_{22}+\delta_{33}$ & Yes \\
    \hline
    T12 & $\tfrac{1}{3}(\delta_{333}+S\delta_{112}+S\delta_{122})$ & Yes & $\delta_{33}+\tfrac{1}{3}(\delta_{11}+\delta_{22}) + \tfrac{4}{3}S\delta_{12}$ & Yes \\
    \hline 
    T13 & $\tfrac{1}{3}(\delta_{222}+S\delta_{113}+S\delta_{133})$ & Yes & $\delta_{22}+\tfrac{1}{3}(\delta_{11}+\delta_{33}) + \tfrac{4}{3}S\delta_{13}$ & Yes \\
    \hline 
    T23 & $\tfrac{1}{3}(\delta_{111}+S\delta_{223}+S\delta_{233})$ & Yes & $\delta_{11}+\tfrac{1}{3}(\delta_{22}+\delta_{33}) + \tfrac{4}{3}S\delta_{23}$ & Yes \\
    \hline
    F112 & $\tfrac{1}{2}S(\delta_{112}+\delta_{233})$ & No & $\tfrac{1}{2}(\delta_{11}+\delta_{33})
+ S(\delta_{12}+\delta_{23})$ & Yes \\
    \hline
    F113 & $\tfrac{1}{2}S(\delta_{113}+\delta_{223})$ & No & $\tfrac{1}{2}(\delta_{11}+\delta_{22})
+ S(\delta_{13}+\delta_{23})$ & Yes \\
    \hline
    F122 & $\tfrac{1}{2}S(\delta_{122}+\delta_{133})$ & No & $\tfrac{1}{2}(\delta_{22}+\delta_{33})
+ S(\delta_{12}+\delta_{13})$ & Yes \\
    \hline
    C & $S \delta_{123}$ & Yes & $S (\delta_{12}+\delta_{13}+\delta_{23})$ & Yes \\
    \hline \hline
    K233,111 & $\tfrac{1}{2}S\delta_{233}+\tfrac{1}{3}\delta_{111}+\tfrac{1}{6}\delta_{222}$ 
& No & $\delta_{11}+\tfrac{1}{2}(\delta_{22}+\delta_{33})
+ S\delta_{23}$ & No \\
    \hline
    K133,222 & $\tfrac{1}{2}S\delta_{133}+\tfrac{1}{3}\delta_{222}+\tfrac{1}{6}\delta_{111}$ 
& No & $\delta_{22}+\tfrac{1}{2}(\delta_{11}+\delta_{33})
+ S\delta_{13}$ & No \\
    \hline
    K122,333 & $\tfrac{1}{2}S\delta_{122}+\tfrac{1}{3}\delta_{333}+\tfrac{1}{6}\delta_{111}$ 
& No & $\delta_{33}+\tfrac{1}{2}(\delta_{11}+\delta_{22})
+ S\delta_{12}$ & No \\
    \hline     
    K223,111 & $\tfrac{1}{2}S\delta_{223}+\tfrac{1}{3}\delta_{111}+\tfrac{1}{6}\delta_{333}$ 
& No & $\delta_{11}+\tfrac{1}{2}(\delta_{22}+\delta_{33})
+ S\delta_{23}$ & No \\
    \hline   
    K113,222 & $\tfrac{1}{2}S\delta_{113}+\tfrac{1}{3}\delta_{222}+\tfrac{1}{6}\delta_{333}$ 
& No & $\delta_{22}+\tfrac{1}{2}(\delta_{11}+\delta_{33})
+ S\delta_{13}$ & No \\
    \hline
    K112,333 & $\tfrac{1}{2}S\delta_{112}+\tfrac{1}{3}\delta_{333}+\tfrac{1}{6}\delta_{222}$ 
& No &
 $\delta_{33}+\tfrac{1}{2}(\delta_{11}+\delta_{22})
+ S\delta_{12}$ & No \\
    \hline
    K112,223 & $S(\tfrac{1}{2}\delta_{112}+\tfrac{1}{4}\delta_{223})+\tfrac{1}{4}\delta_{333}$ 
& No & $\tfrac{3}{4}\delta_{33}+\tfrac{1}{2}\delta_{11}+\tfrac{1}{4}\delta_{22}
+ S(\delta_{12}+\tfrac{1}{2}\delta_{23})$ & No \\
    \hline 
    K113,233 & $S(\tfrac{1}{2}\delta_{113}+\tfrac{1}{4}\delta_{233})+\tfrac{1}{4}\delta_{222}$ 
& No & $\tfrac{3}{4}\delta_{22}+\tfrac{1}{2}\delta_{11}+\tfrac{1}{4}\delta_{33}
+ S(\delta_{13}+\tfrac{1}{2}\delta_{23})$ & No \\
    \hline     
    K223,133 & $S(\tfrac{1}{2}\delta_{223}+\tfrac{1}{4}\delta_{133})+\tfrac{1}{4}\delta_{111}$ 
& No & $\tfrac{3}{4}\delta_{11}+\tfrac{1}{2}\delta_{22}+\tfrac{1}{4}\delta_{33}
+ S(\delta_{23}+\tfrac{1}{2}\delta_{13})$ & No \\
    \hline    
    K122,113 & $S(\tfrac{1}{2}\delta_{122}+\tfrac{1}{4}\delta_{113})+\tfrac{1}{4}\delta_{333}$ 
& No & $\tfrac{3}{4}\delta_{33}+\tfrac{1}{2}\delta_{22}+\tfrac{1}{4}\delta_{11}
+ S(\delta_{12}+\tfrac{1}{2}\delta_{13})$ & No \\
    \hline    
    K133,112 & $S(\tfrac{1}{2}\delta_{133}+\tfrac{1}{4}\delta_{112})+\tfrac{1}{4}\delta_{222}$ 
& No & $\tfrac{3}{4}\delta_{22}+\tfrac{1}{2}\delta_{33}+\tfrac{1}{4}\delta_{11}
+ S(\delta_{13}+\tfrac{1}{2}\delta_{12})$ & No \\
    \hline       
    K233,122 & $S(\tfrac{1}{2}\delta_{233}+\tfrac{1}{4}\delta_{122})+\tfrac{1}{4}\delta_{111}$ 
& No & $\tfrac{3}{4}\delta_{11}+\tfrac{1}{2}\delta_{33}+\tfrac{1}{4}\delta_{22}
+ S(\delta_{23}+\tfrac{1}{2}\delta_{12})$ & No \\
    \hline     
    Id,C & $\tfrac{1}{3}S(\delta_{112}+\delta_{133}+\delta_{223})$ & Yes & $\tfrac{1}{3}(\delta_{11}+\delta_{22}+\delta_{33})
+ \tfrac{2}{3}S(\delta_{12}+\delta_{13}+\delta_{23})$ & No \\
    \hline
    ${\rm Id,C'}$ & $\tfrac{1}{3}S(\delta_{113}+\delta_{122}+\delta_{233})$ & Yes & $\tfrac{1}{3}(\delta_{11}+\delta_{22}+\delta_{33})
+ \tfrac{2}{3}S(\delta_{12}+\delta_{13}+\delta_{23})$ & No \\
    \hline
  \end{tabular}
} 
 \caption{The $22$ extreme points of the symmetric Kantorovich polytope $\calP_{sym,\lambdabar}(X^N)$ for $3$ marginals and $3$ states ($N=\ell=3$). The extreme points fall into two geometric classes, depending on whether their 2-point marginal is also extreme (top $8$ states) or not (bottom $14$ states). Within each class, we have ordered the extreme points by total size (more precisely: sum) of off-diagonal elements of 2-point marginal. Thus elements near the top of their class are expected to optimize attractive costs, while those near the bottom of their class should correspond to repulsive costs. This intuition is made rigorous in Section \ref{S:gallery}.}
\label{Table:1}
\end{table}
We feel that the proof of Theorem \ref{T:1} is not so important and hence we postpone it. Instead we first use the theorem to infer Example 1.1. 
\\[2mm]
{\bf Proof of Example 1.1.} We need to show that the Kantorovich plan $\gamma_*$ or -- in the notation of Theorem \ref{T:1}, Figure \ref{F:KantPoly}, and Table \ref{Table:1} -- $\rmF 112$ is the unique ground state of the optimal transport problem in the example. Recall that in the example,  the competing Kantorovich plans are {\it not} required to be symmetric. We solve the OT problem stepwise, by establishing the following claims: 
\\[1mm]
1) $\gamma_*$ is the unique ground state of $\calC$ on the set $\calP_{sym,\lambdabar}(X^3)$ of symmetric probability measures on $X^3$ with uniform one-point marginal $\lambdabar$. 
\\[1mm]
2) $\gamma_*$ is a ground state of $\calC$ on the set $\calP_{\lambdabar}(X^3)$ of all probability measures on $X^3$ with uniform one-point marginal $\lambdabar$.
\\[1mm]
3) The ground state of $\calC$ on $\calP_{\lambdabar}(X^3)$ is unique.
\\[1mm]
ad 1): By Theorem \ref{T:1} and formula \eqref{costid} it suffices to check that $\calC'[M_2\gamma_*]<\calC'[M_2\tilde{\gamma}]$ where $\tilde{\gamma}$ is any of the extreme points of $\calP_{sym,\lambdabar}(X^3)$ other than $\gamma_*$ whose $2$-point marginal is an extreme point of $\calP_{N-rep,\lambdabar}(X^2)$. From the explicit expressions in Table \ref{Table:1} we can read off the cost in terms of the spring equilibrium bond length $a$ (see Figure \ref{F:FKresults}): 
\begin{eqnarray*}
  \calC[{\rm Id}] &=& 3v(0) \; = \; 3a^2, \\
  \calC[\rmT 12] &=& \calC[\rmT 23] = \tfrac{5}{3} v(0) + \tfrac{4}{3} v(1) \; = \; \tfrac{5}{3}a^2 + \tfrac{4}{3}(1-a)^2, \\
  \calC[\rmT 13] &=& \tfrac{5}{3} v(0) + \tfrac{4}{3} v(2) \; = \; \tfrac{5}{3}a^2 + \tfrac{4}{3}(2-a)^2, \\
  \calC[\rmF 112] &=& v(0) + 2v(1) \; = \; a^2 + 2(1-a)^2, \\
  \calC[\rmF 113] &=& \calC[\rmF 122] \; = \; v(0) + v(1) + v(2) \; = \; a^2 + (1-a)^2 + (2-a)^2, \\
  \calC[\rmC ] &=& 2v(1) + v(2) = 2(1-a)^2 + (2-a)^2.
\end{eqnarray*}  
It follows (see Figure \ref{F:FKresults}) that for $a\in(0,\tfrac{1}{2})$ the unique minimizer is  $\gamma={\rm Id}=\tfrac{1}{3}(\delta_{111}+\delta_{222}+\delta_{333})$, whereas for $a\in(\tfrac{1}{2},1)$, and in particular for $a=\tfrac{3}{4}$, the unique minimizer is $\gamma=\rmF 112 = \tfrac{1}{2}S(\delta_{112}+\delta_{233})$. 
\\[1mm]
ad 2): This follows from the first equality in \eqref{costid}. 
\\[1mm]
ad 3): Let $\gamma$ be any minimizer of $\calC$ on $\calP_{\lambdabar}(X^3)$. By the first equality in \eqref{costid} so is $S\gamma$, whence by 1) we must have $S\gamma=\gamma_*$. Since, for any probability measure on $X^N$, $supp\, \gamma\subseteq supp\, S\gamma$, it follows that $\gamma=a\delta_{112}+b\delta_{121} + c\delta_{211} + d\delta_{233} + e\delta_{323} + f\delta_{332}$ for some $a,\, b, \, c, \, d, \, e, \, f\ge 0$. The marginal condition \eqref{Kantconstr} implies, in particular, that $\gamma(\{ a_1 \}\times X\times X) = \gamma(X\times\{ a_1 \} \times X)=\gamma(X\times X\times \{ a_1 \} ) =\tfrac{1}{3}$, whence $a+b=a+c=b+c=\tfrac{1}{3}$. By elementary linear algebra, this system is uniquely solved by $a=b=c=\tfrac{1}{6}$. Analogously, by evaluating the marginals of $\gamma$ on $\{ a_3 \}$ instead of $\{ a_1 \}$ we obtain $d=e=f=\tfrac{1}{6}$. This shows that $\gamma=\gamma_*$.   
\\[1mm]
Finally we need to show that $\gamma_*$ is neither Monge nor symmetrized Monge. The former is obvious because otherwise $\gamma_*$ would have to be a linear combination of $3$ Dirac measures not $6$. To prove the latter, we note that the relation $\gamma_*=S\gamma$ for some Monge state $\gamma$ would imply, by way of the first equality in \eqref{costid}, that $\gamma$ itself is a minimizer, contradicting uniqueness since $\gamma$ cannot be equal to the non-Monge state $\gamma_*$. 
\\[2mm]
We will infer Theorem \ref{T:1} from the stronger result below which gives the extreme points of the symmetric Kantorovich polytope.
\begin{theorem}\label{T:2} The symmetric Kantorovich polytope $\calP_{sym,\lambdabar}(X^3)$ for $3$ marginals and $3$ states ($N=\ell=3$) has precisely $22$ extreme points, given by the Kantorovich plans in Table \ref{Table:1}. $7$ of these are symmetrized Monge states and $15$ are not (see the 3rd column). 
\end{theorem}
The physical meaning of the $22$ extreme points as molecular packings is shown in Figure \ref{F:KantStatesIconized}. In this visualization, extremality has a simple physical meaning, not limited to $N=\ell=3$:
\be \label{irred}
  \hspace*{-4mm} \begin{array}{c} \mbox{$\gamma$ is an extreme point of the} \\
                   \mbox{symmetric Kantorovich polytope} \end{array}
    \! \Longleftrightarrow 
  \begin{array}{c} \mbox{the molecular packing is {\it irreducible}, that is, not decom-} \\
                   \mbox{posable into uniform fillings with fewer molecules.} \end{array}
\ee
%
%
\noindent
{\bf Proof of Theorem \ref{T:1} using Theorem \ref{T:2}:} Since the reduced polytope $\calP_{N-rep,\lambdabar}(X^2)$ is the image of the symmetric polytope $\calP_{sym,\lambdabar}(X^3)$ under the marginal map $M_2$, and since this map is linear, the reduced polytope is the convex hull of the set
$$
         \{M_2\gamma \, : \, \gamma \mbox{ is an extreme point of }\calP_{sym,\lambdabar}(X^3)\}.
$$ 
The elements of the latter set are listed in Table \ref{T:1} (see column 4). By inspection the first $8$ elements are all different, and by inspection of \ref{F:KantPoly} none of them is contained in the convex hull of the others. Furthermore, we claim that the last $14$ elements in column 4 are strict convex combinations of the first $8$ elements, more precisely:
\begin{align}
 & M_2(Kiij,kkk) = \tfrac{1}{4} M_2({\rm Id}) + \tfrac{3}{4} M_2(Tij) \;\;\; (i, \, j, \, k \mbox{ all different}), \\
 & M_2(Kiij,jjk) = \tfrac{1}{2} \bigl( M_2(\rmF iij) + M_2(\rmK iij,kkk ) \bigr) \;\;\; (i, \, j, \, k \mbox{ all different}), \\
 & M_2({\rm I},\rmC ) = M_2({\rm I}, \rmC ') = \tfrac{1}{3}M_2(I) + \tfrac{2}{3} M_2(\rmC ). 
\end{align}
(Here index sequences like $iij$ are to be interpreted as $jii$ in case $j<i$, to match the notation in Table \ref{Table:1}.) These relations can be seen geometrically in Figure \ref{F:KantPoly2}. 
It follows that the extreme points of the reduced Kantorovich polytope consist of the first $8$ elements of Table \ref{Table:1} column 4, and moreover that these elements have no other inverse image under $M_2$ in the symmetric Kantorovich polytope than the corresponding plans in column 2. This establishes all claims in Theorem \ref{T:1} except which of the first 8 plans in column 2 are symmetrized Monge plans. That the $5$ plans ${\rm Id}$, $\rmT ij$ ($i<j$), and $\rmC$ are symmetrized Monge is clear. The most elementary (if not the most elegant) way to check that $\rmF 122$, $\rmF 113$, $\rmF 122$ are not goes as follows: list the $3!=6$ permutations of $\{1,2,3\}$ in the form $\{\tau_1,...,\tau_6\}$; note that each symmetrized Monge state can be written in the form $S\gamma$ with $\gamma$ as in \eqref{Mongeansatz} with three permutations ${\rm Id}$, $\tau_i$, $\tau_j$ ($i<j$), which gives $\tbinom{6}{2}=15$ possibilities; check that the 10 remaining possibilities other than those leading to the $5$ already discussed Monge plans do not yield $\rmF 122$, $\rmF 113$, $\rmF 122$. A more elegant criterion for when a plan is symmetrized Monge is given in \cite{FV18}. 
\\[2mm]
{\bf Proof of Theorem \ref{T:2}}. The symmetric Kantorovich polytope is a finite-dimensional set defined by both linear inequality and linear equality constraints; by standard results from linear programming (see e.g. \cite{Be09} Section 2) the extreme points can be found by 
solving certain subsystems of the linear constraint equations. In our case these subsystems are at most $3\times 3$ since -- due to the state space consisting only of $3$ points -- the marginal constraint \eqref{Kantconstr} is $3$-dimensional. Thus finding the extreme points is an exercise in linear algebra, and can quickly be done by hand without resorting to (symbolic or numerical) software. For completeness we include the details. 
\\[1mm]
The linear space spanned by the symmetric probability measures on $X^3$ has the canonical basis $\{S\delta_{ijk} \, : \, 1\le i\le j\le k\le \ell\}$, i.e., for $\ell=3$, 
$\{\delta_{111}, \; \delta_{222}, \; \delta_{333}, \; S\delta_{112}, \; S\delta_{113}, \; 
S\delta_{122}, \; S\delta_{223}, $ 
$S\delta_{133}, \; S\delta_{233}, \; S\delta_{123}\}$.  
Thus we may identify symmetric measures on $X^3$ with ten-dimensional vectors
$$
\alpha = (\alpha_{111}, \, \alpha_{222}, \, \alpha_{333}, \, \alpha_{112}, \, \alpha_{113}, \, 
\alpha_{122}, \, \alpha_{223}, \, \alpha_{133}, \, \alpha_{233}, \, \alpha_{123}),
$$
via the expansion $\gamma=\sum_{1\le i\le j \le k}\alpha_{ijk}S\delta_{ijk}$. With this identification, the symmetric Kantorovich polytope \eqref{newset1} turns into the set of vectors $\alpha\in\R^{10}$ satisfying
\begin{align}
    & \;\; \alpha \ge 0 \mbox{ componentwise}, \label{Pos} \\[-10\jot]
    & \underbrace{\begin{pmatrix} 
                  1 & 0 & 0 & \tfrac{2}{3} & \tfrac{2}{3} & \tfrac{1}{3} & 0 & \tfrac{1}{3} & 0 &
                  \tfrac{1}{3} \\[0.5mm]
                  0 & 1 & 0 & \tfrac{1}{3} & 0 & \tfrac{2}{3} & \tfrac{2}{3} & 0 & \tfrac{1}{3} &
                  \tfrac{1}{3} \\[0.5mm]
                  0 & 0 & 1 & 0 & \tfrac{1}{3} & 0 & \tfrac{1}{3} & \tfrac{2}{3} & \tfrac{2}{3} &
                  \tfrac{1}{3} \end{pmatrix}}_{=:A}
    \begin{array}{c} \textcolor{white}{\vdots} \\ \textcolor{white}{.} \\ \begin{pmatrix} \alpha_{111} \\[0.5mm] \alpha_{222} \\[0.5mm] \alpha_{333} \\[0.5mm] \vdots \\[0.5mm] \alpha_{123}  
                     \end{pmatrix} \end{array}
     = \begin{pmatrix} \tfrac{1}{3} \\[0.5mm] \tfrac{1}{3} \\[0.5mm] \tfrac{1}{3} \end{pmatrix}. \label{linconstr}
\end{align} 
It is well known (see e.g. \cite{Be09} Section 2) that the extreme points are given by those vectors satisfying \eqref{Pos}--\eqref{linconstr} such that the matrix columns corresponding to the nonzero components of $\alpha$ are linearly independent. In particular, the number of nonzero components of $\alpha$ must be either $1$, $2$, or $3$, giving us three cases to investigate.
\\[1mm]
{\bf Case 1: $\alpha$ has $1$ nonzero component.} The only constant column of $A$ is the last one, so the nonzero component of $\alpha$ must be the last one, i.e. $\alpha_{123}$. This yields the extreme point $\gamma=S\delta_{123}$ (called $\rmC$ in Table \ref{Table:1} because it is the symmetrized Monge state generated by the powers of the cyclic permutation $\rmC$). 
\\[1mm]
{\bf Case 2: $\alpha$ has $2$ nonzero components.} First we claim that $\alpha_{123}$ must be zero. This is because otherwise the constraint \eqref{linconstr} would force the column of $A$ corresponding to the other nonzero component of $\alpha$ to be a constant vector, contradicting linear independence. Next, the first $3$ components of $\alpha$ must also be zero, because otherwise $A\alpha$ would contain two unequal components. Thus the columns associated with the nonzero components of $\alpha$ must both come from the ``middle block'' (columns 4 to 9). The only combinations that work are $\alpha_{112}=\alpha_{233}=\tfrac{1}{2}$, $\alpha_{113}=\alpha_{223}=\tfrac{1}{2}$, and $\alpha_{122}=\alpha_{233}=\tfrac{1}{2}$, yielding the three ``Frenkel-Kontorova'' extreme points $\rmF 112$, $\rmF 113$, $\rmF 122$ in the table. 
\\[1mm]
{\bf Case 3: $\alpha$ has $3$ nonzero components.} First, note that again $\alpha_{123}$ must be zero, since otherwise the marginal condition forces the sum of the other two columns to be a constant vector, contradicting linear independence. We subdivide the remaining possibilities according to the number of columns from the ``first block'' (columns 1 to 3). 
\\[1mm]
{\bf 3a: $3$ columns from the first block.} In this case $\alpha_{111}=\alpha_{222}=\alpha_{333}=\tfrac{1}{3}$, yielding $\gamma=\tfrac{1}{3}(\delta_{111}+\delta_{222}+\delta_{333})$ (called ${\rm Id}$ in the table because it is the Monge state generated by the powers of the identity). 
\\[1mm]
{\bf 3b: $2$ columns from the first block.} The remaining column must be one of the six columns from the middle block. If it is, e.g., the first of these columns (corresponding to $\alpha_{112}\neq 0$), one finds that the other two columns must be the ones corresponding to $\alpha_{222}$ and $\alpha_{333}$, yielding $\gamma=\tfrac{1}{2}S\delta_{112} + \tfrac{1}{6}\delta_{222} + \tfrac{1}{3}\delta_{333}$. Proceeding analogously for the other five columns from the middle block gives the six states from $\rmK 233,111$ to $\rmK 112,333$ in the table. 
\\[1mm]
{\bf 3c: $1$ column from the first block.} Suppse, e.g., that it is the first column, corresponding to $\alpha_{111}\neq 0$. The other two columns cannot contain the entry $\tfrac{2}{3}$, because otherwise the first component of $A\alpha$ would exceed the maximum of the other two components. This leaves only $\tbinom{4}{2}=6$ choices. The only ones yielding a solution to \eqref{Pos}--\eqref{linconstr} are the column pairs corresponding to $\alpha_{122}$ and $\alpha_{233}$, or $\alpha_{133}$ and $\alpha_{223}$, or $\alpha_{223}$ and $\alpha_{233}$. This gives the three extreme points 
$\tfrac{1}{4}(\delta_{111}+S\delta_{122}) + \tfrac{1}{2}S\delta_{233} = \rmK 233,122$, 
$\tfrac{1}{4}(\delta_{111}+S\delta_{133}) + \tfrac{1}{2}S\delta_{223} = \rmK 223,133$, and
$\tfrac{1}{3}(\delta_{111}+S\delta_{223}+S\delta_{233}) =\rmT 23$. Proceeding analogously for $\alpha_{222}$ and $\alpha_{333}$ yields the six states from $\rmK 112,223$ to $\rmK 233,122$ and the three states from $\rmT 12$ to $\rmT 23$ in the table. 
\\[1mm]
{\bf 3d: No column from the first block.} Hence all three columns come from the middle block. Suppose, e.g., that one of them is the first column, corresponding to $\alpha_{112}\neq 0$. This leaves only $\tbinom{5}{2}=10$ choices for the additional column pair. The only choice yielding a solution to \eqref{Pos}--\eqref{linconstr} is the pair corresonding to $\alpha_{223}$ and $\alpha_{113}$, giving $\gamma=\tfrac{1}{3}S(\delta_{112} + \delta_{223} + \delta_{133})$. Taking into account that we could have started with any other component from the middle block in place of $\alpha_{112}$, we obtain the two states ${\rm I,C}$ and ${\rm I,C'}$ in the table. This completes the derivation of the set of extreme points.
\\[1mm]
Finally, to complete the proof of Theorem \ref{T:2} it remains to check which of these extreme points is symmetrized Monge. This can be done in a straightforward manner by proceeding as in the proof of Theorem \ref{T:1}.

\section{Convex geometry of the Monge ansatz} 
The Monge ansatz \eqref{Mongeansatz} corresponds to replacing the reduced Kantorovich polytope of $N$-representable two-point probability measures with uniform marginal by the convex hull of those probability measures which come from symmetrized Monge states,
\be \label{Mongeset}
  \calM_{N-rep}(X^2) = conv\{p\in\calP(X^2)\, : \, p=M_2S\gamma \mbox{ for some }\gamma \mbox{ of form }\eqref{Mongeansatz}\},
\ee
where $conv$ means convex hull. We call this set the {\it reduced Monge polytope for $N$ marginals and $\ell$ states}. Unfortunately, the definition \eqref{Mongeset}, just like \eqref{Plambdabar}, is very abstract, and poorly understood. Many questions pose themselves:
\\[2mm]
1) What does the passage from \eqref{Plambdabar} to \eqref{Mongeset} mean {\it geometrically}? 2) How ``big'' is the difference between the two polytopes? 3) How big is the induced error in optimal cost, at worst? 4) For precisely which class of costs is the error zero? 
\\[2mm]
In our model case $N=\ell=3$, the reduced Monge polytope is straightforwared to calculate explicitly (see Proposition \ref{P:Monge} below), and is plotted in Figure \ref{F:MongePoly}, thereby answering Question 1). One sees that as compared to the Kantorovich polytope in Figure \ref{F:KantPoly}, three slices of significant size around each of the vertices $\rmF 112$, $\rmF 113$, $\rmF 122$ are cut off. 
\begin{figure}
\begin{center}
\includegraphics[width=0.75\textwidth]{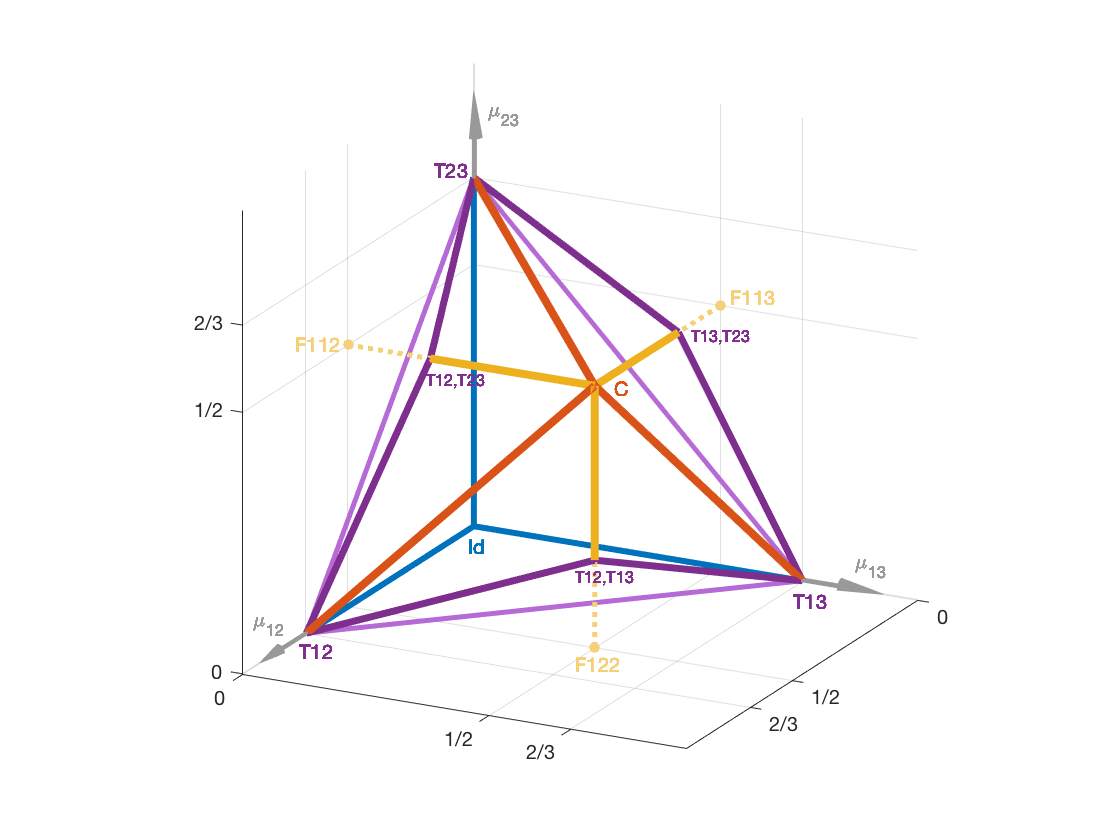}
\end{center}
\caption{The reduced Monge polytope for $3$ marginals and $3$ states ($N=\ell=3$). It is obtained from the Kantorovich polytope in Figure \ref{F:KantPoly} by moving the non-Monge vertices $\rmF 112$, $\rmF 113$, $\rmF 122$ by an amount of $\tfrac{1}{6}$ towards the point $C$.}
\label{F:MongePoly}
\end{figure}
\begin{proposition} \label{P:Monge} The reduced Monge polytope \eqref{Mongeset} has precisely $8$ extreme points, given by $M_2\gamma$ where $\gamma$ is one of the following plans: 
\begin{align}
 & {\rm Id, \; T12, \; T13, \; T23, \; C} \;\; \mbox{(see Table \ref{Table:1})}, \\
 & \underbrace{S\tfrac{1}{3}(\delta_{121}+\delta_{213} + \delta_{332})}_{=:{\small \rmT 12,\!\rmT 23}},  \;\;\; 
   \underbrace{S\tfrac{1}{3}(\delta_{131}+\delta_{223} + \delta_{312})}_{=:{\small \rmT 13,\rmT 23}}, \;\;\;
   \underbrace{S\tfrac{1}{3}(\delta_{123}+\delta_{212} + \delta_{331})}_{=:{\small \rmT 12,\rmT 13}}. 
\end{align}
The last three states correspond, respectively, to the collection of permutations $\{{\rm Id},\, \rmT 12, \, \rmT 23\}$, $\{{\rm Id},\, \rmT 13, \, \rmT 23\}$, $\{{\rm Id},\, \rmT 12, \, \rmT 13\}$ (see Remark 3 in the introduction for nomenclature). 
\end{proposition}

\begin{figure}
\begin{center}
\includegraphics[width=0.6\textwidth]{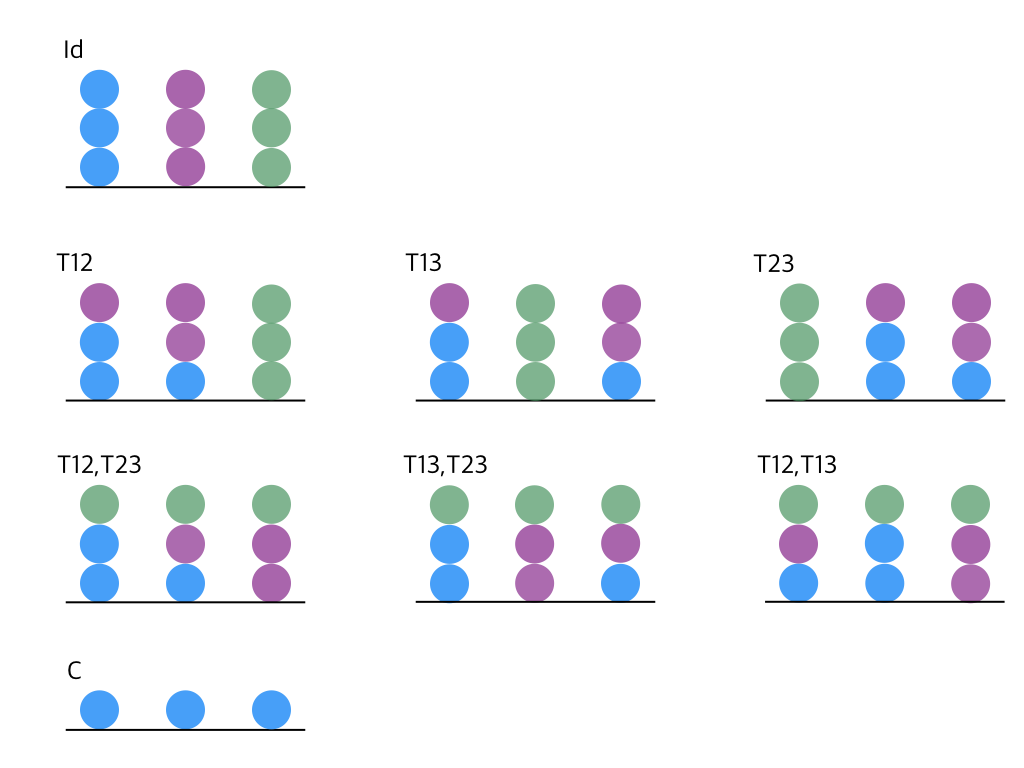}
\end{center}
\caption{Physical meaning of the $8$ extreme points of the reduced Monge polytope $\calM_{N-rep}(X^2)$ for $3$ marginals and $3$ states ($N=\ell=3$) as molecular packings. This picture immediately reveals that the states in the 3rd row fail to be extreme points of the symmetric Kantorovich polytope as they can be decomposed into two uniform packings with fewer molecules. This fact does not appear to be so obvious from the mathematical expressions in Proposition \ref{P:Monge}.}
\label{F:MongeStatesIconized}
\end{figure}

{\bf Proof} The proof is straightforward: as in the proof of Theorem \ref{T:1}, we write each symmetrized Monge state in the form $S\gamma$ with $\gamma$ as in \eqref{Mongeansatz} and three permutations ${\rm Id}$, $\tau_i$, $\tau_j$ ($i<j$), which gives 15 possibilities, compute the associated $2$-point marginals, and check which of them lie outside the convex hull of the others. 
\\[2mm]
The physical meaning of the extreme points is shown in Figure \ref{F:MongeStatesIconized}. 
\section{Gallery of 3-marginal problems on 3 sites with pairwise interaction minimized by Monge states} \label{S:gallery}
Our geometric results from Section \ref{S:geo} can not just be used to find counterexamples to the Monge ansatz, but also to justify it for various interesting costs. 
Let $X=\{a_1,a_2,a_3\}$, i.e. $\ell=3$, and let $d$ be an arbitrary metric on $X$. Prototypical is the euclidean distance $d(x,y)=|x-y|$ when $X$ is a collection of points in $\R^d$. Our interest is in the three-marginal OT problem with pairwise costs which depend only on the metric distances between the points in $X$ (i.e., physically, the interparticle distances). That is to say we consider functionals of form
\be \label{metrcost}
   \calC[\gamma ] = \int_{X^3} \bigl( v(d(x,y)) + v(d(y,z)) + v(d(x,z)) \bigr) \, d\gamma(x,y,z) \mbox{ for some potential }v\, : \, [0,\infty )\to\R.
\ee
Our goal is to minimize \eqref{metrcost} over the set \eqref{newset1} of symmetric probability measures on $X^3$ with uniform marginal, for various interesting classes of potentials $v$. Of course, it would be of interest to extend the results below to arbitrary finite state spaces \eqref{space}, i.e. an arbitrary number of sites, or even to continuous state spaces. Nevertheless we feel that analyzing the minimalistic case of only 3 sites is of value because one can see without being impeded by combinatoric or analytic complications how common consitutive assumptions on the potential such as monotonicity or convexity drive the minimizers towards the ``Monge'' region of the symmetric Kantorovich polytope.    
\\[2mm]
{\bf Example 4.1: Positive power cost $v(d)=d^p$, $0<p<\infty$}. Note that for $p=2$ this is the Gangbo-\'{S}wi\k{e}ch cost \cite{GS98}. The unique minimizer of \eqref{metrcost} on $\calP_{sym,\lambdabar}(X^3)$ is the Monge state $\gamma=\tfrac{1}{3}(\delta_{111} + \delta_{222} + \delta_{333})$, which corresponds to the collection of permutations $\{\tau_1,\tau_2,\tau_3\} = \{{\rm Id},{\rm Id},{\rm Id}\}$. This is trivial because the integrand in \eqref{metrcost} is minimized pointwise if and only if $x=y=z$, and the above $\gamma$ is the only element of $\calP_{sym,\lambdabar}(X^3)$ which is supported on the diagonal $\{(x,x,x)\, : \, x\in X\}$.   
\\[2mm]
{\bf Example 4.2: Negative power cost $v(d)=\tfrac{1}{d^{\alpha}}$, $\alpha>0$}. Note that for $\alpha=1$ this is the Coulomb cost \cite{CFK11,BDG12,FMPCK13,Pa13,CDD13,CFM14,CFP15}, see also \cite{SGS07,Se99}. The unique minimizer of \eqref{metrcost} on $\calP_{sym,\lambdabar}(X^3)$ 
is the symmetrized Monge state $\gamma=S\delta_{123}$. This is because $v(d(x,x))=\infty$ for all $x\in X$ (since $d$ is a metric, which implies $d(x,x)=0$), and the above $\gamma$ is the unique element of $\calP_{sym,\lambdabar}(X^3)$ whose 2-point marginal vanishes on the diagonal $\{(x,x)\, : \, x\in X\}$ of $X^2$. 
\\[2mm]
{\bf Example 4.3: General attractive costs, i.e. $v(d)$ strictly increasing in $d$}. As in Example 4.1, the unique minimizer is the Monge state $\gamma=\tfrac{1}{3}(\delta_{111} + \delta_{222} + \delta_{333})$ generated by powers of the identity. This follows by arguing as in Example 4.1.
\\[2mm]
{\bf Example 4.4: General repulsive costs, i.e. $v(d)$ strictly decreasing in $d$}. This example is a little less trivial, and illustrates the virtues of the Monge ansatz. We claim that there always exists a minimizer of Monge form. This can be seen as follows. Denote $d_{ij}:=d(a_i,a_j)$, $c_{ij}:=v(d_{ij})$, $\mu_{ij}:=\ell (M_2\gamma)(\{(a_i,a_j)\} )$. For any $\gamma\in \calP_{sym,\lambdabar}(X^3)$ we find by first using the symmetry of $\gamma$ and then the marginal equations $\mu_{11}=1-(\mu_{12}+\mu_{13})$, 
$\mu_{22}=1-(\mu_{12}+\mu_{23})$, $\mu_{33}=1-(\mu_{13}+\mu_{23})$ that the cost \eqref{metrcost} evaluates to 
\begin{eqnarray}
  \calC[\gamma] &=& \sum_{i=1}^3 c_{ii}\mu_{ii} + 2\sum_{i<j} c_{ij}\mu_{ij} \; = \; 
   c_{11} + c_{22} + c_{33} \\
   && + \bigl(2c_{12}-(c_{11}\! +\! c_{22})\bigr) \mu_{12} 
                               + \bigl(2c_{13}-(c_{11}\! +\! c_{33})\bigr) \mu_{13} 
                               + \bigl(2c_{23}-(c_{22}\! +\! c_{33})\bigr) \mu_{23}.  
\label{Cred}
\end{eqnarray}
Since, for $i\neq j$, $d_{ij}>0=d_{ii}=d_{jj}$, and $v$ is strictly decreasing, we have 
\be \label{star}
    2c_{ij} - (c_{ii}+c_{ij}) < 0 \mbox{ for all }i<j.
\ee
This means that the vertex $(\mu_{12},\mu_{23},\mu_{13})=(\tfrac12,\tfrac12,\tfrac12)$ of the reduced Kantorovich polytope (see Figure \ref{F:KantPoly}) (corresponding to $\gamma=\rmC$) always gives lower cost than the vertex $(0,0,0)$ (corresponding to $\gamma={\rm Id}$) and the vertices 
  $(\tfrac12,\tfrac12,0)$, $(\tfrac12,0,\tfrac12)$, $(0,\tfrac12,\tfrac12)$ (corresponding to $\rmF 112$, $\rmF 122$, $\rmF 113$). This shows that for any strictly decreasing $v$, one of $\rmT 12$, $\rmT 13$, $\rmT 23$, $\rmC$ is a minimizer. In particular, there always exists a minimizer of symmetrized Monge form, as asserted. Moreover every minimizer is a convex combination of these 4 states. That all these 4 states can indeed occur as unique minimizers of repulsive costs can be seen from Examples 4.5 and 4.6. 
\\[2mm]
{\bf Example 4.5: Repulsive positive-power cost, i.e. $v(d)=-d^p$, $p>0$}. Here we confine ourselves to $X=\{1,2,3\}\subset\R$ (equidistant points on the line), $d(x,y)=|x-y|$ (euclidean 1D metric). Note that $p=2$ corresponds to the repulsive harmonic oscillator \cite{DGN15, GKR18}. Heuristically, because the largest distance is $d_{13}$, if $p$ is large this cost should favour the vertex $(\mu_{12},\mu_{23},\mu_{13})=(0,0,\tfrac{2}{3})$ of the reduced Kantorovich polytope, i.e. the state $\rmT 13$. Indeed, we already know from Example 4.4 that any minimizer is a convex combination of the four states $\rmT 12$, $\rmT 13$, $\rmT 23$, $\rmC$, and by plugging them into \eqref{Cred} we find that for $p<\log 6/\log 2=2.58496...$ the unique minimizer of \eqref{metrcost} on $\calP_{sym,\lambdabar}(X^3)$ is $\rmC$, while for $p>\log 6/\log 2$ it is $\rmT 13$. The critical exponent emerges because $\calC[\rmC]-\calC[\rmT 12]=-2 + 2^p/3$. 
\\[2mm]
{\bf Example 4.6: General repulsive convex costs, i.e. $v$ strictly decreasing and convex}. 
Here we consider again an arbitrary metric $d$ on $X$. We claim that the unique minimizer of \eqref{metrcost} on $\calP_{sym,\lambdabar}(X^3)$ is the symmetrized Monge state $\gamma=S\delta_{123}$ generated by the cyclic permutation $\rmC$. By the findings of Example 4.4, it suffices to show that $\calC[\rmC]<\calC[\rmT ij]$ for all $i<j$. Consider, for instance, $\rmT 13$. To start with, by the triangle inequality for the metric, 
$$
   d_{13}\le d_{12} + d_{23} =: d'_{13},
$$
and consequently -- since $v$ is decreasing -- 
\be \label{one}
   c_{13} = v(d_{13}) \ge v(d_{13}')=:c_{13}'.
\ee
Moreover, since $d_{13}'>0$ we may write, with $\lambda := d_{12}/d_{13}'$, $d_{12}=(1-\lambda)\cdot 0 + \lambda d_{13}'$ and $d_{23}=\lambda\cdot 0 + (1-\lambda)d_{13}'$. By convexity of $v$, it follows that
\begin{eqnarray*}
   c_{12} &\le & (1-\lambda )v(0) + \lambda c_{13}', \\
   c_{23} &\le & \lambda v(0) + (1-\lambda )c_{13}'.
\end{eqnarray*}
Adding these inequalities yields 
\be \label{two}
   c_{12} + c_{23} \le v(0) + c_{13}'.
\ee
Armed with these inequalities, we can now analyze the cost difference $\calC[\rmC]-\calC[\rmT 13]$. By \eqref{Cred} and the fact that $\rmC$ and $\rmT 13$ correspond, respectively, to $(\mu_{12},\mu_{23},\mu_{13})=(\tfrac12,\tfrac12,\tfrac12)$ and $(0,0,\tfrac23)$ (see Figure \ref{F:KantPoly}), we find that 
$$
   \calC[\rmC]-\calC[\rmT 13] = (2c_{12}-2v(0))\cdot \tfrac12 + (2c_{23}-2v(0))\cdot\tfrac12 + 
         (2c_{13}-2v(0))\cdot (-\tfrac16) = c_{12}+c_{23}-\tfrac13 c_{13} - \tfrac53 v(0).
$$
Plugging \eqref{one}, \eqref{two} into this expression gives
$$
  \calC[\rmC] - \calC[\rmT 13] \le \tfrac{2}{3}(c_{13}-v(0)),
$$
which is $<0$ since $v$ is strictly decreasing. Arguing analogously for $\rmT 12$ and $\rmT 23$ establishes that $\rmC$ is the unique minimizer.
\section{Frenkel-Kontorova gas and formation of microstructure}
\label{S:micro}
\begin{figure}
\begin{minipage}{0.67\textwidth}
\begin{center}
\includegraphics[width=\textwidth]{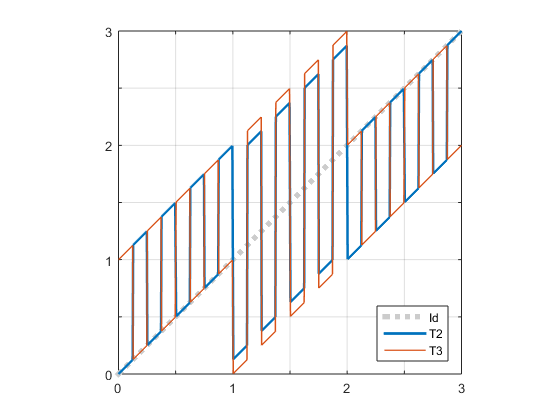} 
\end{center}
\end{minipage}
\begin{minipage}{0.32\textwidth}
\begin{center}
\includegraphics[width=0.76\textwidth]{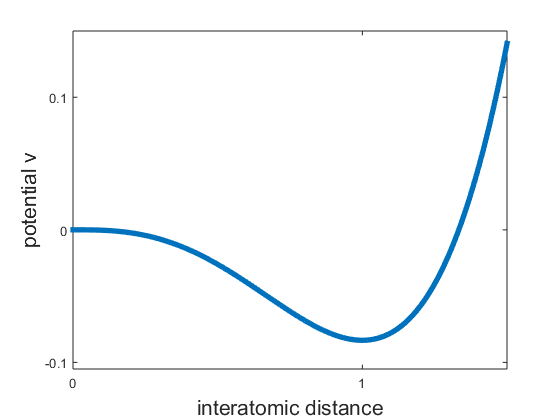} 
\end{center}
\end{minipage}
\vspace*{-5mm}
\caption{Microstructure and nonattainment in Monge multi-marginal optimal transport. {\it Left panel}: Minimizing sequence of the problem in Example \ref{E:cts}. Physically the problem corresponds to finding the ground state of a 1D homogeneous gas of 3-particle molecules, with particles (``atoms'') within a molecule interacting via the potential $v$. {\it Right panel}: The  Frenkel-Kontorova-like interatomic potential $v(r)=\tfrac{r^4}{4}-\frac{r^3}{3}$ which causes these  wild oscillations.
}
\label{F:micro}
\end{figure}

At first sight our finite-state-space results in Section \ref{S:geo} may seem far from continuous optimal transport. In fact they can be used to establish the following simple -- and physically intuitive -- continuous counterexample to the Monge ansatz in 3-marginal optimal transport. 
\\[2mm]
Let us replace the finite state space \eqref{space} by an interval $X=[a,b]$ and the marginal measure \eqref{unif} by its continuous analogue on $[a,b]$, 
\be \label{unif'}
    \lambda = \tfrac{1}{b-a} \, dx \; \mbox{ where $dx$ is the Lebesgue measure}. 
\ee
For $3$ marginals ($N=3$), the marginal condition \eqref{Kantconstr} becomes
\be \label{Kantconstr'}
    \gamma( \Omega\times X\times X) = \gamma( X\times\Omega\times X) = \gamma( X\times X\times\Omega) \mbox{ for all Borel sets }\Omega\subseteq X
\ee
(notation: $\gamma\mapsto\lambda$), leading to the set of Kantorovich plans $\calP_\lambda(X^3)=\{\gamma \in \calP(X^3)\, : \, \eqref{Kantconstr'}\}$. Here $\calP(X^3)$ denotes the set of Borel probability measures on $X^3$. 
The Monge ansatz \eqref{Mongeansatz} turns into
\be \label{Monge'}
    \gamma = \int_{[a,b]} d(x) \delta_{T_2(x)}(y) \delta_{T_3(x)}(z) \, dx \mbox{ for some }T_2, \, T_3\in\calM_\lambda
\ee
(or, in optimal transport notation, $\gamma = (id,T_2,T_3)_\sharp \lambda$),  
where $d(x)=\tfrac{1}{b-a}$ is the (constant) density of the measure \eqref{unif'} and $\calM_\lambda$ is the usual set of Monge maps which push $\lambda$ forward to itself,  
\be \label{MongeMaps}
    \calM_\lambda = \{ T \, : \, X\to X \, : \, T \mbox{ Borel measurable, }T_\sharp \lambda = \lambda \}. 
\ee
Here, as usual, $(T_i)_\sharp \lambda = \nu$ means that $T$ pushes $\lambda$ forward to $\nu$, i.e. $\lambda(T^{-1}(\Omega)) = \nu(\Omega)$ for all Borel sets $\Omega$. Maps belonging to \eqref{MongeMaps} are called Monge maps and the resulting 3-point probability measures \eqref{Monge'} are called Monge plans. As before our interest is in symmetric and pairwise costs
\be \label{Kantcts}
  \calC[\gamma ]= \int_{[a,b]^3} \bigl( v(|x-y|) + v(|y-z|) + v(|x-z|) \bigr) \, d\gamma(x,y,z), \;\; v\, : \, [0,\infty )\to\R \mbox{ continuous}.
\ee
For Monge plans \eqref{Monge'} the cost obviously satisfies $\calC[\gamma]=I[T_2,T_3]$ where 
\be \label{Mongered}
   I[T_2,T_3] = \int_a^b d(x) \bigl( v(|x-T_2(x)|) + v(|T_2(x)-T_3(x)|) + v(|x-T_3(x)|) \bigr)
   \, dx .
\ee
\begin{example} \label{E:cts} (1D homogeneous Frenkel-Kontorova gas) Let $X=[0,3]$ and
$v(r)=\frac{r^4}4-\frac{r^3}3$, and consider the 3-marginal OT problem of minimizing \eqref{Kantcts} subject to $\gamma\mapsto\lambda$. (Physically, this corresponds to seeking the ground state of a 1D homogeneous gas of 3-particle molecules, with particles within a molecule interacting via the Frenkel-Kontorova-like potential $v$ (see Figure \ref{F:micro}, right panel.) \\[1mm]
a) The infimum of $\calC$ over the set of Monge plans \eqref{Monge'} is not attained. The infimum over symmetrized Monge plans is not attained either. \\[1mm]  
b) A minimizer of $\calC$ over the set of Kantorovich plans $\calP_\lambda(X^3)$ is given by 
$$
   \gamma = S\int_0^3 \!\! \alpha(x) \delta_{T_2(x)}(y)\delta_{T_3(x)}(z) \, dx
$$
(or, in optimal transport notation, $\gamma = S \, (id,T_2,T_3)_\sharp  \, \alpha \, dx$) where $S$ is the symmetrization operator \eqref{S} and  
$$
   \alpha(x) = \begin{cases} \tfrac12, & x\in [0,1]\cup[2,3] \\
                               0,      & x\in (1,2), \end{cases}, \;\; 
   T_2(x)=x, \;\; T_3(x) = \begin{cases} x+1, & x\in[0,1], \\ 
                                        x,   & x\in(1,2) \\
                                      x-1, & x\in[2,3]. \end{cases}
$$
\end{example}
We remark that the plan in b) has the form of a continuous analogue of a {\it quasi-Monge state}; the latter notion is introduced in \cite{FV18}. 
\\[2mm]
To understand heuristically what is going on, let us rewrite the Kantorovich minimizer as 
\be \label{super}
   \gamma = \int_{[0,1]} S \Bigl( \tfrac12 \delta_a (x) \delta_a (y) \delta_{a+1}(z) 
                                 +\tfrac12 \delta_{a+1} (x) \delta_{a+2} (y) \delta_{a+2}(z) \Bigr) \, da. 
\ee
This is nothing but a continuous superposition of translates of the Frenkel-Kontorova minimizer $\gamma_*=\rmF 112$ from Example \ref{E:counter}. The first contribution in the integrand is a copy of the blue molecule in the right panel of Figure \ref{F:FKmodel}, translated by $a$ to the right; the second contribution is a copy of the purple molecule, translated likewise.   
\\[2mm]
As regards the nonattainment of the Monge infimum, the reader may wonder what a typical minimizing sequence, i.e. a sequence of Monge maps $T_2^{(\nu)}$, $T_3^{(\nu)}$ such that $\lim_{\nu\to\infty} I[T_2^{(\nu)},T_3^{(\nu)}] = \inf I$, will look like. 
An example of a minimizing sequence is shown in Figure \ref{F:micro}. (To obtain the sequence from the single picture, take the width of the intervals where the maps are affine to be $\tfrac{1}{2\nu}$.) 
Note that if the width of the intervals on which the maps are affine is made smaller and smaller, the triple of distances $|x-T_2(x)|$, $|T_2(x)-T_3(x)|$, $|x-T_3(x)|$ approaches the optimal set of values $\{(0,1,1),(1,0,1),(1,1,0)\}$, for every $x$. But this means that $T_2^{(\nu)}$ and $T_3^{(\nu)}$ develop faster and faster oscillations (``microstructure'') everywhere in the domain  and converge weakly but not strongly in $L^p([0,3])$ ($1<p<\infty$) to the non-minimizing maps 
$$
   T_2(x)=T_3(x)=\begin{cases} x+\tfrac12, & x\in[0,1] \\
                               x,          & x\in(1,2) \\
                               x-\tfrac12, & x\in[1,2]. \end{cases}
$$             
Unfortunately we cannot offer an analogous result for {\it all} minimizing sequences. 
\\[2mm]
{\bf Proof of the assertions in Example \ref{E:cts}}. The proof is perhaps even more instructive than the example itself. Instead of ad hoc arguments we proceed by reduction to a finite state space, then use the geometric results from Section \ref{S:geo}. 

First of all let us determine the pointwise minimizers of the integrand in \eqref{Kantcts}. Because $x,y,z\in\R$, one of the distances $|x-y|$, $|y-z|$, $|x-z|$ is the sum of the other two, and the elementary calculus problem of minimizing $v(r_1)+v(r_2)+v(r_1+r_2)$ on $[0,\infty)^2$ is solved precisely when $(r_1,r_2)=(0,1)$ or $(1,0)$ (hence our choice of the potential $v$). It follows that the integrand in \eqref{Kantcts} is minimized if and only if $(|x-y|,|y-z|,|x-z|)\in\{ (0,1,1),(1,0,1),(1,1,0)\}$. The plan \eqref{super} is clearly supported on this set, establishing b). 
To show a), note first that by b), any minimizer $\tilde{\gamma}$ must be supported in the above set. Suppose $\tilde{\gamma}$ is of Monge form. Since $|T_j(x)-x|\in\{0,1\}$ a.e., we must in particular have $T_j(x)\in[0,3]\cap (x+\Z)$ a.e. We partition the integration in \eqref{Mongered} 
according to these finitely many possiblities, and parametrize the possiblities by two maps $t_1,\, t_2\, : \, B=\{-1,0,1\} \to B$. Why this is a good parametrization will be come clear shortly. For any two such maps, define
$$
     \Omega_{t_2,t_3} := \{ x\in (1,2) \, : \, T_2(x+b)=x+t_2(b) \mbox{ for all }b\in B, \; 
                                               T_3(x+b)=x+t_3(b) \mbox{ for all }b\in B\} .
$$ 
Then, up to a set of measure zero, 
$$
    [0,3] = \bigcup_{t_2, \, t_3\, : \, B\to B}  \Omega_{t_2,t_3} \cup (\Omega_{t_2,t_3}-1) \cup (\Omega_{t_2,t_3}+1) \;\; (\mbox{disjoint union})
$$
and consequently, by decomposing the integral over $[0,3]$ in \eqref{Mongered} accordingly,   
\be \label{partition}
   \calC[\tilde{\gamma} ] = \sum_{t_2,t_3\, : \, B\to B} |\Omega_{t_2,t_3}| \, \tilde{\calC}[\tilde{\gamma}_{t_2,t_3}], \;\;\; \sum_{t_2,t_3\, : \, B\to B} |\Omega_{t_2,t_3}|=1,
\ee
where $\tilde{\gamma}_{t_2,t_3}$ is the following probability measure on the finite state space $B^3$
$$
     \tilde{\gamma}_{t_2,t_3} = \tfrac13 \sum_{b\in B=\{-1,0,1\} } \delta_{b} \otimes \delta_{t_2(b)} \otimes \delta_{t_3(b)},
$$
and $\tilde{C}$ is the finite-state-space cost on $\calP(B^3)$ obtained from \eqref{Kantcts} by replacing the domain of integration $[a,b]^3$ by $B^3$. In fact we know more about the $\tilde{\gamma}_{t_2,t_3}$. By the push-forward condition in \eqref{MongeMaps} we have 
$(T_j)_\sharp \lambda|_{\Omega\cup(\Omega-1)\cup(\Omega+1)} = \lambda|_{\Omega\cup(\Omega-1)\cup(\Omega+1)}$ for all Borel $\Omega\subseteq (1,2)$ and consequently, whenever $|\Omega_{t_2,t_3}|>0$, $(t_j)_\sharp \lambdabar = \lambdabar$ where $\lambdabar$ is the uniform measure $\tfrac{1}{3}\sum_{b\in B}\delta_b$ on $B$, that is to say the $t_j$ are permutations. So the $\tilde{\gamma}_{t_2,t_3}$ are finite Monge states; but we know from Section \ref{S:geo} that the cost of symmetrized Monge states (and hence, by \eqref{minid}, of Monge states) is strictly bigger than the Kantorovich minimum of $\tilde{\calC}$, since the only extreme point of the symmetric Kantorovich polytope supported in the set of minimizers of $v(|x-y|)+v(|y-z|)+v(|x-z|)$ on $B^3$ is the non-Monge plan $\rmF 112$ (see Theorem \ref{T:2} or Figure \ref{F:KantStatesIconized}). This together with \eqref{partition} shows that the minimum cost in the continuous and the finite problem are equal, i.e.
$$
   \min_{\gamma\in\calP_\lambda([0,3]^3)} \calC[\gamma] = \min_{\tilde{\gamma}\in\calP_\lambdabar(B^3)} \tilde{\calC}[\tilde{\gamma}],
$$ 
and that $\tilde{\gamma}$ as well as its symmetrization $S\tilde{\gamma}$ have strictly bigger cost, i.e.
\be \label{ineq}
   \calC[\tilde{\gamma}] = \calC[ S\tilde{\gamma} ] > \min_{\gamma\in\calP_\lambda ([a,b]^3)} \calC[\gamma].
\ee
Finally, it is well known for much more general OT problems than the one here that the Kantorovich minimum equals the Monge infimum. Together with \eqref{ineq} this establishes b). 
\\[5mm]
{\bf Acknowledgements.} The author thanks Maximilian Fichtl and Daniela V\"ogler for helpful discussions. 

\begin{small}

\end{small}

\end{document}